\documentclass[11pt]{article}
\usepackage[numbers,sort&compress]{natbib}
\usepackage{enumerate}
\usepackage{amscd}
\usepackage{amsmath}
\usepackage{latexsym}
\usepackage{amsfonts}
\usepackage{setspace}
\usepackage{amssymb}
\usepackage{amsthm}
\usepackage{verbatim}
\usepackage{mathrsfs}
\usepackage{enumerate}
\usepackage[hypertexnames=false]{hyperref}
\usepackage[utf8]{inputenc}
\usepackage{amsmath, amssymb, pgfplots}
\pgfplotsset{compat = newest, width = 12cm, height =12cm}
\usepackage{tikz}
\usepackage{caption}

\theoremstyle{plain}
\theoremstyle{definition}\newtheorem{theorem}{Theorem}[section]
\theoremstyle{definition}\newtheorem{lemma}[theorem]{Lemma}
\theoremstyle{plain}
\theoremstyle{plain}\newtheorem{prop}[theorem]{Proposition}
\theoremstyle{definition}\newtheorem{remark}{Remark}[section]
\usepackage{xcolor}
\def\v{\varepsilon}

\newcommand{\e}{\mathrm{e}}

\allowdisplaybreaks

\numberwithin{equation}{section}
\begin{document}
	%%%%%%%%%%%%%%%%%%%%%%%%%%%%%%%%%%%%%%%%%%%%%%%%%%%%%%%%%%%%%%%%%%%%%%%%%%%%%%%%%%%%%%%%%%%%%%%%%%%%
	\title{Periodic Solutions to the Steady Viscous Burgers Equation: A Constructive Example}
	\author{Quansen Jiu\footnote{School of Mathematical Sciences, Capital Normal University, Beijing, 100048, P. R. China. Email: jiuqs@cnu.edu.cn}~~~\,\,\,\,Yuhan Cao\footnote{School of Mathematical Sciences, Capital Normal University, Beijing, 100048, P. R. China. Email: caoyuhan@163.com}}
	\date{}
	\maketitle
	\begin{abstract}
		In this paper, we consider the periodic problem to steady viscous Burgers equation with an periodic and external force. Our main aim is to construct periodic solutions to this model which are uniformly bounded with respect to the viscosity via a direct Fourier series approach. Firstly, as an example, starting from a special solution of the non-viscous Burgers equation with a specific force, we solve  the approximate solutions to the viscous Burgers equation in an explicit way. Secondly,  we construct the solution to the periodic problem of the viscous Burgers equation with the external force by solving the initial value problem to a second ordinary differential equations, which is uniformly bounded with respect to the viscosity. Compared with Jauslin-Kreiss-Moser’s result in \cite{M}, our approach  provides an explicit constructive procedure for the periodic solutions and yields a  sharper convergence rate (up to higher order) when the viscosity vanishes.
		\end{abstract}

	\tableofcontents
	
	%%%%%%%%%%
	\section{Introduction and main results}
	\qquad	The one-dimensional viscous Burgers equation can be written as
	\begin{equation}\label{01}
		u_t+uu_x-\varepsilon u_{xx}=0,
	\end{equation}
	where $u=u(x,t)$ is an unknown function, $\varepsilon \geq0$ represents the viscosity coefficient. The Burgers equation is a well-known hyperbolic-parabolic equation and characterizes a nonlinear dissipative system. It is used as a simplified model to describe many physical phenomena.
	The study of the Burgers equation has a long history starting with the seminal
	papers by Burgers \cite{B}, Cole \cite{C} and Hopf \cite{H} in which the Cole-Hopf transformation was introduced. When $\varepsilon=0,$ equation \eqref{01} becomes the classical non-viscous Burgers equation
	\begin{equation}\label{02}
		u_t+uu_x=0.
	\end{equation}\\
	The non-viscous equation \eqref{02} is a bases for study of non-linear waves in PDE (see \cite{S}).
	
	More recently there have been several articles dealing with the forced Burgers(see \cite{K-L} and references therein)
	equation
	\begin{equation}\label{03}
		u_t-\varepsilon u_{xx}+uu_x=f,
	\end{equation}
	where $f$ is a given function of $x$ and $t$ which represents an external force. Over the past two decades, the Burgers equation\eqref{03} with time- and space-periodic solutions (denoted by $\mathrm{Z}^{2}$-periodic solutions if the periodic in space and time are both one) has been investigated. Among the notable contributions are the works of Jauslin, Kreiss, and Moser who showed the existence and uniqueness of the solution to \eqref{03} for some specific smooth forcing function $ f(t,x)$(see \cite{M}). More precisely, it is proved in \cite{M} that
		\begin{prop}\label{prop1}
		Assume that $f=F_x(t,x)$, where $F(t,x)$ is a smooth function of period one in both variables $t$ and $x$. Then equation \eqref{03} with $\varepsilon >0$ has a unique $\mathrm{Z}^{2}$-periodic solution satisfying
		$\int_0^1 u ~dx=c$ for any fixed $c \in \mathrm{R}$.
	\end{prop}

	When the external force is time-independent, the equation \eqref{03} admits a  solution satisfying the steady viscous Burgers equation,
	\begin{equation}\label{04}
		\frac{1}{2}(U^2)_x=\varepsilon U_{xx}+F_x.
	\end{equation}
	Moreover, the solution is uniformly bounded on the viscosity $\varepsilon$. More precisely, it is shown in \cite{M} that
	\begin{prop}\label{prop2}
		Assume that $f=F_x(x)$ depends on $x$ only. Then for any $\varepsilon>0$, the solution presented in Proposition \ref{prop1} converges, for $ t \rightarrow \infty$, to a unique 1-periodic solution $U(x,\varepsilon):=u_{\infty}^{(\varepsilon)}(x)$ of the steady viscous Burgers equation \eqref{04} satisfying
		$$\int_0^1U dx=c.$$
		Morever, there is a constant $M>0$ which does not depend on $\varepsilon$ such that
		\begin{equation*}
			|U| \leq M , U_x \leq M,  \int_0^1|U_x| dx \leq M.
		\end{equation*}
	\end{prop}

	It is noted that Jauslin, Kreiss, and Moser also proved that there exists the unique eternal trajectory(attractor) which attracts in future all solutions of the problem.
	In \cite{E}, Weinan E addressed the open problem from \cite{M}, providing a further characterization of the viscous and periodic solutions to the forced Burgers equation \eqref{03}.
    The work of \cite{M} was later extended in \cite{Chen}\cite{V}\cite{Z}, among them, \cite{Z} eliminated the requirement for the time periodicity of the forcing term and further proves that for any $L^{\infty}((-\infty,\infty);L^2(0,1))$ forcing, a unique eternal trajectory attracts all trajectories. For more references please reffer to \cite{FTKS}\cite{Cy}\cite{CZ} which are concerned with computer-assisted researches on equation \eqref{03} and \eqref{04}.
	
	In this paper, we consider the following one-dimensional steady viscous Burgers equation:
	\begin{equation}\label{1.1}
		uu_{x}-\varepsilon u_{xx}=f(x), ~x \in [0,2\pi]
	\end{equation}
	where $\varepsilon>0$ is the viscosity coefficient, $u=u(x)$  is the unknown function,
	$f(x)$ is a given $2\pi$-periodic function with respect to $x$.
	Periodic conditions to \eqref{1.1} are imposed as $$u(0)=u(2\pi).$$
	The aim of this paper is to construct a periodic solution of one-dimensional steady Burgers equation\eqref{1.1} which is uniformly bounded with respect to the viscosity $\varepsilon$.
	Denote the solution of equation \eqref{1.1} by $u^{\varepsilon}(x)$, which is depend on $\varepsilon$ in general.
	Then we will construct the solution to \eqref{1.1} with the following formal expansion
	\begin{equation}\label{1.2}
		u^\v (x)=u_0(x)+\v u_1(x)+\cdots+\v ^nu_n(x)+\cdots .
	\end{equation}
	Substituting \eqref{1.2} into \eqref{1.1}, we can obtain recurrence relations between the approximate solutions $u_n(x) (n=1,2,\cdots)$, which is
	\begin{equation*}
		u_0(x)u_k(x)=u_{k-1}^\prime(x)-\frac{1}{2}\sum_{i+j=k}u_i(x)u_j(x) , ~k=1,2,\cdots.
	\end{equation*}
	
	Given some special starting function $u_0(x)$ and $f(x)$, we can solve out the approximate solutions explicitly.
	In the first part of this paper, as an example, we take
	$$u_0(x)=-(2+\cos x), ~f=u_0(x)u_{0x}=-2\sin x-\sin x\cos x.$$
	Through Fourier series expansion, we can give explicit expressions of $u_1$ and $u_2$ as follows.
	$$u_1(x)=\sum\limits_{n=1}^\infty 2(-2+\sqrt{3})^n\sin nx.$$  $$u_2(x)=\sum\limits_{n=0}^\infty -(-2+\sqrt{3})^n\left(\frac{\sqrt{3}}{2}n^2+\frac{n}{3}+\frac{\sqrt{3}}{18}\right) \cos nx. $$
	For the cases $n=3,4,\cdots$, the expressions of the approximate solutions are given as follows .
	\begin{theorem}\label{0-1}		Let $m$ be any positive integer. Assume that for any $k \in \mathrm{N}_{+}$ and $k\leq m-1$, one has
		\begin{equation}
			\begin{split}\label{0.1}
				u_k=&\frac{1+(-1)^k}{2}[\sum_{n=0}^{\infty}t^n(C_{k,2k-2}n^{2k-2}+C_{k,2k-3}n^{2k-3}+\cdots +C_{k,0})\cos nx]\\
				+&\frac{1-(-1)^k}{2}[\sum_{n=0}^{\infty}t^n(C_{k,2k-2}n^{2k-2}+C_{k,2k-3}n^{2k-3}+\cdots +C_{k,0})\sin nx],
			\end{split}
		\end{equation}
		where $C_{k,j}$ are some coefficients with $j=0,1,2,\cdots,2k-2$.\\
		Then when $k=m,$ the expression \eqref{0.1} is also ture.
	\end{theorem}
	It is remarked that Theorem \ref{0-1} was first proved in \cite{J}. For the sake of completeness and self-containment, we will give a sketch of proof of Theorem \ref{0-1} in this paper.
	
	In the second part of this paper, we prove that for any positive integer $n \in \mathrm{N}_{+}$, the solution to equation \eqref{1.1} admits the following expansion
	$$u^{\v}(x)=u_0(x)+\v u_1(x)+\cdots+\v^n u_n(x)+o(\v^n).$$
	More precisely, we get the following theorem
	\begin{theorem}\label{th-11.2}
		There exisis $\v_0>0$, which depends on $n$, such that for any $0<\v<\v_0$, the equation \eqref{1.1} has a periodic solution $u^\v(x) \in C^2([0,2\pi])$, which satisfies
		$$\|u^\v(x)-u_0(x)-\v u_1(x)-\cdots-\v^nu_n(x)\| \leq C\v^{n+1},$$
	\end{theorem}
	where $C>0$ is a constant depending on $n$ but independent of $\v$.
	
	To prove Theorem \ref{th-11.2},
	we introduce an error function $v_{n+1}(\frac{x}{\v})$, satisfying
	$$u^{\v}(x)=u_0(x)+\v u_1(x)+\cdots+\v^n u_n(x)+\v^{n+1}v_{n+1}(\frac{x}{\v}).$$
	Then we get the error equation
	\begin{equation}\label{a}
		\begin{aligned}
			&-v_{n+1}^{\prime\prime}(\frac{x}{\v})+\v^{n+1} v_{n+1}(\frac{x}{\v})v_{n+1}^{\prime}(\frac{x}{\v})+
			\big(u_0(x)+\v u_1(x)+\v^2 u_2(x)+\dots+\v^n u_n(x)\big)v_{n+1}^{\prime}(\frac{x}{\v})+\\
			&\big(\v u_0^\prime(x)+\v^2u_1^\prime(x)+\v^3 u_2^{\prime}(x)+\dots+\v^{n+1} u_n^{\prime}(x)\big) v_{n+1}(\frac{x}{\v})+\\
			&\sum_{k=1}^{n}{\v^k(u_k(x)u_n^\prime(x)+u_{k+1}(x)u_{n-1}^\prime(x)+\dots+u_n(x)u_k^\prime(x))}-\v u_n^{\prime\prime}(x)=0.
		\end{aligned}
	\end{equation}
	Let
	$$y=\frac{x}{\v},$$
	it deduces that
	\begin{equation}\label{b}
		\begin{aligned}
			&-v_{n+1}^{\prime\prime}(y)+\v^{n+1} v_{n+1}(y)v_{n+1}^{\prime}(y)
			+\big(u_0(\v y)+\v u_1(\v y)+\dots+\v^n u_n(\v y)\big)v_{n+1}^{\prime}(y)\\
			&+\big(\v u_0^\prime(\v y)+\v^2u_1^\prime(\v y)+\v^3 u_2^{\prime}(\v y)+\dots+\v^{n+1} u_n^{\prime}(\v y)\big) v_{n+1}(y)
			+F(\v y)=0.
		\end{aligned}
	\end{equation}
	To prove Theorem \ref{th-11.2}, we only need to prove
	\begin{theorem}\label{th-1.3}
		There exists $\v_{0}>0$ such that for any $0< \v <\v_{0}$, the equation \eqref{b}
		has a solution $v_{n+1}(y) \in C^2([0,\frac{2\pi}{\v}])$, satisfying $$v_{n+1}(0)=v_{n+1}(\frac{2\pi}{\v})$$
		and
		$$|v_{n+1}| \leq C,$$
		where $C>0$ is a constant depending on $n$ but independent of $\v$.
	\end{theorem}	
	\begin{remark}
		Compared with Proposition \ref{prop2}, we present an explicit and constructive approximate solutions to the periodic solutions of \eqref{1.1}.  Moreover, we derive a sharper convergence rate up to any order through Fourier expansions.
	\end{remark}
	We now present the main idea of the proof of Theorem \ref{th-1.3}. We first construct the solution to the initial value problem of \eqref{b}, that is for any given initial data $v_{n+1}(0)=x_{n+1},
	v_{n+1}^{\prime}(0)=y_{n+1},$ we consider the Cauchy problem
	\begin{equation}\label{c}
		\left\{
		\begin{array}{l}
			\begin{aligned}
				&-v_{n+1}^{\prime\prime}(y)+\v^{n+1} v_{n+1}(y)v_{n+1}^{\prime}(y)+\big(u_0(\v y)+\v u_1(\v y)+\dots+\v^n u_n(\v y)\big)v_{n+1}^{\prime}(y)\\
				&+\big(\v u_0^\prime(\v y)+\v^2u_1^\prime(\v y)+\v^3 u_2^{\prime}(\v y)+\dots+\v^{n+1} u_n^{\prime}(\v y)\big) v_{n+1}(y)
				+F(\v y)=0,\\
				&v_{n+1}(0)=x_{n+1},\\
				&v_{n+1}^{\prime}(0)=y_{n+1}.
			\end{aligned}	\\
		\end{array}
		\right.
	\end{equation}
	It is easy to obtain that the Cauchy problem \eqref{c} has a unique solution $v_{n+1}(x) \in [0,T_{n+1})$, where $T_{n+1}>0$ is a constant denpending on $\v$(maybe small).
	
	Secondly, we prove that the existence interval of solution to \eqref{c} is at least $[0,\frac{2\pi}{\v}]$, by using upper and lower solutions method. Meanwhile, we obtain that the solution of \eqref{c} is uniformly bounded.
	
	Finally,  based on the solution of \eqref{c} on $[0,\frac{2\pi}{\v}]$, we construct the periodic solution of \eqref{b} satisfying $v_{n+1}(0)=v_{n+1}(\frac{2\pi}{\v})$. To this end, we denote also the solution of \eqref{c} by $v_{n+1}(y,y_{n+1})$. A key observation is that the periodic condition is equivalent to
	$$\Phi_{n+1}(\v,y_{n+1}):=-v_{n+1}^{\prime}(\frac{2\pi}{\v},y_{n+1})+\frac{\v^{n+1}}{2}v_{n+1}^2(\frac{2\pi}{\v},y_{n+1})+y_{n+1}-\frac{\v^{n+1}}{2}x_{n+1}^2=0.$$
	Our goal is to apply the implicit function theorem to determine the relationship between $y_{n+1}$ and $\v$ such that $\Phi_{n+1}(\v, y_{n+1})=0.$ To achieve this, we need to extend the definition of $\Phi$ to include the points where $\v=0$. We introduce $r_{n+1}(y,y_{n+1})=v_{n+1}(y,y_{n+1})-\tilde{v}_{n+1}(y,y_{n+1})$, where $\tilde{v}_{n+1}(y,y_{n+1})$ satisfies the non-viscous equation
	\begin{equation}\label{e2}
		\left\{
		\begin{array}{ll}
			-\tilde{v}_{n+1}^{\prime}(y,y_{n+1})-3\tilde{v}_{n+1}(y,y_{n+1})+3x_{n+1}+y_{n+1}=0,\\
			\tilde{v}_{n+1}(0,y_{n+1})=x_{n+1}.
		\end{array}
		\right.		
	\end{equation}
	\eqref{e2} can be solved out to get
	$$\displaystyle{\tilde{v}_{n+1}(y,y_{n+1})=x_{n+1}+\frac{y_{n+1}}{3}(1-e^{-3y}).}$$
	Then $r_{n+1}(y,y_{n+1})$ satisfies the following equation
	\begin{equation}\label{e}
		\left\{
		\begin{array}{ll}
			-r_{n+1}^{\prime}(y,y_{n+1})-3r_{n+1}(y,y_{n+1})+\frac{\v^{n+1}}{2} v_{n+1}^2(y,y_{n+1})+\big(3+u_0(\v y)+\dots+\v^{n} u_n(\v y)\big)v_{n+1}\\
			+F_1(\l^2 y)-\frac{\v^{n+1}}{2}x_{n+1}^2-\big(\v u_1(0)+\v^2 u_2(0)+\dots+\v^{n} u_n(0)\big)x_{n+1}=0,\\
			r_{n+1}(0)=0.
		\end{array}
		\right.		
	\end{equation}
	Since $v_{n+1}$ is uniformly bounded in $\v$, we will prove that
	$r_{n+1}(\v)$ converges to $0$ as $\v \rightarrow 0$ uniformly with respect to $y_{n+1}$. Then we can define
	$$	\Phi_{n+1}(0,y_{n+1})=\lim\limits_{\v \rightarrow 0}\Phi_{n+1}(\v,y_{n+1})=y_{n+1}.$$
	Then $\Phi_{n+1}(0,y_{n+1})$ is continuous and $\Phi_{n+1}(0,0)=0.$
	Similarly, we obtain $\frac{\partial \Phi_{n+1}}{\partial y_{n+1}}$ is continuous and $\frac{\partial \Phi_{n+1}}{\partial y_{n+1}}(0,0)=1 \neq 0.$ Then by the implicit function theorem, we get the existence of periodic solutions of error equation \eqref{a}.
	
	The paper is organized as follows. In Section 2, we construct an approximate solution by giving recurrence relations. In Section 3, we  establish the existence of periodic solutions for the first-order and higher-order error equations, which leads to the proof of our main results.

	\section{Construction of the approximate solutions $\mathbf{u_n}$ }
	In this section, we will derive a recurrence relation between the approximate solutions $u_n(n=1,2,\cdots)$, and give a specific constructive method of $u_n$.
	
	Differentiating both sides of \eqref{1.2} with respect to $x$ yields
	\begin{align*}
		u^\varepsilon _x(x)&=u_{0x}(x)+\varepsilon u_{1x}(x)+\cdots+\varepsilon ^nu_{nx}(x)+\cdots, \\
		u^\varepsilon _{xx}(x)&=u_{0xx}(x)+\varepsilon u_{1xx}(x)+\cdots+\varepsilon ^nu_{nxx}(x)+\cdots.
	\end{align*}
	Substituting them into the \eqref{1.1}, we obtain
	\begin{equation}\label{2.1}
		\begin{aligned}
			&\varepsilon[u_0(x)u_{1x}(x)+u_{0x}(x)u_1(x)-u_{0xx}(x)]\\
			+&\varepsilon^2[u_0(x)u_{2x}(x)+u_1(x)u_{1x}(x)+u_2(x)u_{0x}(x)-u_{1xx}(x)]\\
			+&\cdots\\
			+&\varepsilon^k[u_0(x)u_{kx}(x)+u_1(x)u_{(k-1)x}(x)+\cdots+u_{k-1}(x)u_{1x}(x)+u_k(x)u_{0x}(x)-u_{(k-1)xx}(x)]\\
			+&\cdots=0 .
		\end{aligned}
	\end{equation}
	It follows that
	\begin{equation}\label{22.2}
		\left\{
		\begin{array}{ll}
			&u_0(x)u_{1x}(x)+u_{0x}(x)u_1(x)-u_{0xx}(x)=0 , \\
			&u_0(x)u_{2x}(x)+u_1(x)u_{1x}(x)+u_2(x)u_{0x}(x)-u_{1xx}(x)=0 , \\
			&\vdots\\
			&u_0(x)u_{kx}(x)+u_1(x)u_{(k-1)x}(x)+\cdots+u_k(x)u_{0x}(x)-u_{(k-1)xx}(x)=0, \\
			&\vdots\\
		\end{array}
		\right.
	\end{equation}
	Integrating each equation in \eqref{22.2} with respect to $x$  yields (Note that the constant term is taken to be $0$ during integration).
	\begin{equation}\label{22.3}
		\left\{
		\begin{array}{ll}
			&u_0(x)u_1(x)=u_0^\prime(x), \\
			&u_0(x)u_2(x)=u_1^\prime(x)-\frac{1}{2}u^2_1(x), \\
			&\vdots\\
			&u_0(x)u_{2k}(x)=u^\prime_{2k-1}(x)-[u_1(x)u_{2k-1}(x)+\cdots+u_{k-1}(x)u_{k+1}(x)+\frac{1}{2}u^2_k(x)], \\
			&u_0(x)u_{2k+1}(x)=u^\prime_{2k}(x)-[u_1(x)u_{2k}(x)+u_2(x)u_{2k-1}(x)+\cdots+u_k(x)u_{k+1}(x)], \\
			&\vdots\\
		\end{array}
		\right.
	\end{equation}
	It follows from \eqref{22.3} that for any positive integers $i$, $j$, and $k$,
	\begin{equation}\label{2.2}
		u_0(x)u_k(x)=u_{k-1}^\prime(x)-\frac{1}{2}\sum_{i+j=k}u_i(x)u_j(x) .
	\end{equation}
	It is easy to verify that if $u_0(x)$ is an even function, then $u_{2k+1}(x)$ is an odd function and $u_{2k}(x)$ is an even function for $k=1,2,\cdots.$
	
	As an example, we take $u_0(x)=-(2+\cos x),~f(x)=u_0(x)u_{0x}(x)=-\sin x\cos x-2\sin x$.
	We will first calculate $u_1, u_2, u_3$ respectively.
	
	\subsection{The expression of $\mathbf{u_1}$}
	Since $u_1(x)$ is an odd function, it can be expressed as
	$$u_1(x)=b_{1,1}\sin x +b_{1,2}\sin 2x +\cdots+b_{1,n}\sin nx+\cdots .$$
	The first equation in \eqref{22.3} becomes
	\begin{equation}\label{2.3}
		-(\cos x +2)(b_{1,1}\sin x +b_{1,2}\sin 2x +\cdots+b_{1,n}\sin nx+\cdots)=\sin x .
	\end{equation}
	Using
	\begin{equation*}
		\sin \alpha \cos \beta =\frac{\sin (\alpha+\beta)+\sin (\alpha-\beta)}{2} ,
	\end{equation*}
	comparing the coefficients on both sides of \eqref{2.3}, we obtain
	\begin{equation}\label{2.4}
		\left\{
		\begin{array}{ll}
			2b_{1,1}+\frac{1}{2}b_{1,2}=-1 , \\
			\frac{1}{2}b_{1,1}+2b_{1,2}+\frac{1}{2}b_{1,3}=0 , \\
			\frac{1}{2}b_{1,2}+2b_{1,3}+\frac{1}{2}b_{1,4}=0 , \\
			\vdots\\
			\frac{1}{2}b_{1,n-1}+2b_{1,n}+\frac{1}{2}b_{1,n+1}=0 ,\\
			\vdots\\
		\end{array}
		\right.
	\end{equation}
	Namely, for positive integers $n\geq2$, we have
	\begin{equation}\label{2.5}
		b_{1,n+1}+4b_{1,n}+b_{1,n-1}=0 .
	\end{equation}
	Rewrite \eqref{2.5} as
	$$b_{1,n+1}-[(-2+\sqrt{3})+(-2-\sqrt{3})]b_{1,n}+[(-2+\sqrt{3})(-2-\sqrt{3})]b_{1,n-1}=0 . $$
	Let $t=-2+\sqrt{3}, q=-2-\sqrt{3}$, and rewrite the above equation to obtain $$b_{1,n+1}-tb_{1,n}=q[b_{1,n}-tb_{1,n-1}],~n\geq2. $$
	Set
	\begin{equation}\label{2.6}
		b_{1,2}-tb_{1,1}=0 .
	\end{equation}
		Then, based the first equation of \eqref{2.4}, we obtain
	$$b_{1,1}=2t,~b_{1,n}=2t^n,~n=2, 3, \cdots.$$
	Thus we have $u_1(x)=\sum\limits_{n=1}^\infty 2t^n\sin nx $ which
	converges absolutely and uniformly on $[0,2\pi]$.\\
	\begin{remark}
		It should be remarked that \eqref{2.6} is a necessary and sufficient condition for the convergence of the Fourier series of $u_1(x)$.
		In fact, if \eqref{2.6} dose not hold, suppose $b_{1,1}=2t+\varepsilon(\varepsilon\neq0)$, denote the perturbation term by $\tilde{b}_{1,1}=\varepsilon$, then by the first equation of \eqref{2.4}, we have $\tilde{b}_{1,2}=-4\varepsilon.$
		Hence, $\{\tilde{b}_{1,n}\}$ satisfies the system of equations
		\begin{equation*}
			\left\{
			\begin{array}{ll}
				\tilde{b}_{1,n+1}+4\tilde{b}_{1,n}+\tilde{b}_{1,n-1}=0 , \\
				\tilde{b}_{1,1}=\varepsilon ,~~~~~~ \tilde{b}_{1,2}=-4\varepsilon . \\
			\end{array}
			\right.
		\end{equation*}
		Solving the above system of equations, we obtain
		$$\tilde{b}_{1,n}=\frac{t^n-q^n}{2\sqrt{3}}\cdot\varepsilon, ~~n=2,3,\cdots. $$
		Note that $\sum\limits_{n=1}^\infty \frac{t^n}{2\sqrt{3}}$ is convergent, but $\sum\limits_{n=1}^\infty \frac{q^n}{2\sqrt{3}}$ is divergent, so $\sum\limits_{n=1}^\infty \tilde{b}_{1,n}$ is divergent. Since $b_{1,n}=2t^n+\tilde{b}_{1,n}$, it deduces that $\sum\limits_{n=1}^\infty b_{1,n}$ is also divergent.
	\end{remark}

	\subsection{The expression of $\mathbf{u_2}$}
	By \eqref{2.2}, $u_2$ satisfies
	\begin{equation}\label{2.7}
		u_0u_2=u_1^\prime-\frac{1}{2}u_1^2 .
	\end{equation}
	Note that the Fourier expansion of $u_2$ is
	\begin{equation*}
		u_2(x)=\frac{a_{2,0}}{2}+a_{2,1}\cos x+a_{2,2}\cos 2x+\cdots+a_{2,n}\cos nx+\cdots.
	\end{equation*}
	Substituting it into \eqref{2.7} yields\
	\begin{equation}\label{2.8}
		-(\cos x+2)(\frac{a_{2,0}}{2}+a_{2,1}\cos x+\cdots+a_{2n}\cos nx+\cdots)=\sum_{n=1}^{\infty}nb_{1,n}\cos nx-\frac{1}{2}(\sum_{n=1}^{\infty}2t^n\sin nx)^2 .
	\end{equation}
	
	The coefficients on the left-hand side of \eqref{2.8} satisfy
	\begin{itemize}
		\item The constant term is $\displaystyle{}splaystyle{-a_{2,0}-\frac{a_{2,1}}{2}};$
		\item the coefficient of $\cos nx $ is $\displaystyle{}splaystyle{-2a_{2,n}-\frac{a_{2,n-1}+a_{2,n+1}}{2}}  $ for $n\geq1 $.
	\end{itemize}
	The coefficients on the right-hand side of \eqref{2.8} satisfy
	\begin{itemize}
		\item The constant term is $\displaystyle{}splaystyle{-\frac{t^2}{1-t^2}} ;$
		\item the coefficient of $\cos nx $ is $\displaystyle{}splaystyle{t^n(-\frac{2t^2}{1-t^2}+3n-1)} $ for $n\geq1 $.
	\end{itemize}
	Let $A_{2,n}=2t^n(-\frac{2t^2}{1-t^2}+3n-1)$. Using \eqref{2.8}, we have
	\begin{equation}
		\left\{
		\begin{array}{ll}\label{2.9}
			-2a_{2,0}-a_{2,1}=-\frac{2t^2}{1-t^2} , \\
			-4a_{2,n}-a_{2,n-1}-a_{2,n+1}=A_{2,n}, n\geq1.
		\end{array}
		\right.
	\end{equation}
	In view of the second equation of \eqref{2.9}, we introduce  $d_{2,n}$ with $n\geq1$ such that
	\begin{align}
		\begin{cases}
			-a_{2,n+1}+ta_{2,n}-d_{2,n}=q(-a_{2,n}+ta_{2,n-1}-d_{2,n-1}) , \\
			d_{2,n}-qd_{2,n-1}=A_{2,n}.\label{2.10}
		\end{cases}
	\end{align}
	By the second equation of \eqref{2.10}, it deduces
	\begin{equation*}
		\begin{split}
			d_{2,n}&=q^n\left(d_{2,0}+\frac{A_{2,1}}{q}+\cdots+\frac{A_{2,n}}{q^n}\right),~n=1,2,\cdots.
		\end{split}
	\end{equation*}
	Set \begin{equation}\label{2.11}
		d_{2,0}=-\sum_{k=1}^\infty \frac{A_{2,k}}{q^k} .
	\end{equation}
	Then it follows that
	\begin{equation}
		\begin{split}
			d_{2,n}&=q^n\left(\frac{-A_{2,n+1}}{q^{n+1}}+\frac{-A_{2,n+2}}{q^{n+2}}+\cdots\right)\\
			&=\frac{-A_{2,n+1}}{q}+\frac{-A_{2,n+2}}{q^2}+\cdots+\frac{-A_{2,n+k}}{q^k}+\cdots\\
			&=t^{n+1}\left(\sqrt{3}n+\frac{3\sqrt{3}+2}{6}\right),~n=1,2,\cdots. \label{22.14}
		\end{split}
	\end{equation}
	Note that it follows from \eqref{2.11} that $d_{2,0}=\frac{5-4\sqrt{3}}{6}$ which is consistent with \eqref{22.14}.
	\\Now we solve
	\begin{equation*}
		\begin{cases}
			-2a_{2,0}-a_{2,1}&=-\frac{2t^2}{1-t^2} , \\
			-a_{2,1}+ta_{2,0}&=d_{2,0} ,
		\end{cases}
	\end{equation*}
	and obtain $$ a_{2,0}=-\frac{\sqrt{3}}{18}.$$
	In \eqref{2.10}, we choose  $a_{2,n}=ta_{2,n-1}-d_{2,n-1},n\geq1$, it follows that
	\begin{equation*}
		\begin{split}
			a_{2,n}&=ta_{2,n-1}-d_{2,n-1}\\
			&=t^na_{2,0}-t^{n-1}d_{2,0}-\cdots-d_{2,n-1}\\
			&=t^{n}a_{2,0}-t^n\left(\frac{\sqrt{3}}{2}n^2+\frac{n}{3}\right)\\
			&=-t^n\left(\frac{\sqrt{3}}{2}n^2+\frac{n}{3}+\frac{\sqrt{3}}{18}\right), n\geq0.
		\end{split}
	\end{equation*}
	It is easy to prove that $\sum\limits_{n=0}^\infty a_{2,n}\cos nx$ is absolutely and uniformly convergent on $[0,2\pi]$.\\
	
	\begin{remark}
		It should be remarked that \eqref{2.11} is a necessary and sufficient condition for the convergence of the Fourier series of $u_2(x)$ as well. In fact, if \eqref{2.11} does not hold, suppose $d_{2,0}=-\sum\limits_{k=1}^\infty \frac{A_{2,k}}{q^k}+\varepsilon(\varepsilon\neq0)$, denote the perturbation term by $\tilde{d}_{2,0}=\varepsilon$.
		Combine the first equation of \eqref{2.9} with the first equation of \eqref{2.10} when $n=1$, we obtain
		\begin{align}\nonumber
			\begin{cases}
				2\tilde{a}_{2,0}+\tilde{a}_{2,1}&=0, \\
				\tilde{a}_{2,1}-t\tilde{a}_{2,0}&=-\tilde{d}_{2,0}.
			\end{cases}
		\end{align}
		Solve the equations, we have
		\begin{align}\nonumber
			\begin{cases}
				\tilde{a}_{2,0}=\frac{\varepsilon}{2+t}, \\
				\tilde{a}_{2,1}=-\frac{2\varepsilon}{2+t}.
			\end{cases}
		\end{align}
		Then by \eqref{2.10}, we get $\tilde{d}_{2,n}=q^{n}\varepsilon.$
		Since $a_{2,n}=ta_{2,n-1}-d_{2,n-1},n\geq1$, iterating yields
		\begin{align*}
			\tilde{a}_{2,n}&=-\frac{2t^{n-1}}{2+t} \varepsilon-qt^{n-2}\varepsilon-q^{2}t^{n-3}\varepsilon-\cdots-q^{n-2}t\varepsilon-q^{n-1}\varepsilon \\
			&=-q\frac{q^{n-1}-t^{n-1}}{q-t} \varepsilon-\frac{2t^{n-1}}{2+t} \varepsilon.
		\end{align*}
		Noting that $\sum\limits_{n=1}^\infty \frac{t^n}{2\sqrt{3}}$ is convergent, but $\sum\limits_{n=1}^\infty \frac{q^n}{2\sqrt{3}}$ is divergent, then $\sum\limits_{n=0}^\infty \tilde{a}_{2,n}$ is divergent. That is to say, \eqref{2.11} is the only choice that ensures the convergence of the Fourier series of $u_{2}(x)$.
	\end{remark}

	\subsection{The expression of $\mathbf{u_3}$}
	By \eqref{2.2}, $u_3$ satisfies
	\begin{equation}\label{2.12}
		u_0u_3=u'_2-u_1u_2 .
	\end{equation}
	Since $u_{3}$ is odd, we assume
	\begin{equation}\label{2.13}
		u_3(x)=\sum_{n=1}^{\infty}b_{3,n}\sin nx .
	\end{equation}
	For convenience, we express $u_1$ and $u_2$ as
	\begin{align*}
		u_1(x)&=\sum_{n=0}^{\infty}t^nC_{1,0}\sin nx, \\
		u_2(x)&=\sum_{n=0}^{\infty}t^n(C_{2,2}n^2+C_{2,1}n+C_{2,0})\cos nx.
	\end{align*}
	Substituting them into \eqref{2.12} gives
	\begin{equation}
		\begin{split}\label{2.14}
			&-(\cos x+2)(b_{3,1}\sin x+b_{3,2}\sin 2x+b_{3,3}\sin 3x+\cdots)\\
			=&-\sum_{n=0}^{\infty}t^n(C_{2,2}n^3+\cdots+C_{2,0}n)\sin nx-(\sum_{n=0}^{\infty}t^nC_{1,0}\sin nx)[\sum_{n=0}^{\infty}t^n(C_{2,2}n^2+\cdots+C_{2,0})\cos nx].
		\end{split}
	\end{equation}
	The coefficients on the left-hand side of \eqref{2.14} satisfy
	\begin{itemize}
		\item the coefficient of $\sin x $ is $\displaystyle{}splaystyle{-2b_{3,1}-\frac{b_{3,2}}{2}} ; $
		\item the coefficient of $\sin nx$ is $\displaystyle{}splaystyle{-2b_{3,n}-\frac{b_{3,n+1}+b_{3,n-1}}{2}}$ for   $n\geq2 . $
	\end{itemize}
	For the right-hand side of \eqref{2.14}, employing a similar method to calculate $u_2$, we can conclude that the coefficients on the right-hand side of \eqref{2.14} can be assumed as $t^n(D_{3,3}n^3+D_{3,2}n^2+D_{3,1}n+D_{3,0})$.
	Let
	$A_{3,n}:=-2t^n(D_{3,3}n^3+D_{3,2}n^2+D_{3,1}n+D_{3,0})$. Then we have
	\begin{equation}\label{2.15}
		4b_{3,1}+b_{3,2}=A_{3,1} ,
	\end{equation}
	\begin{equation}\label{2.16}
		b_{3,n+1}+4b_{3,n}+b_{3,n-1}=A_{3,n},\ (n\geq2).
	\end{equation}
	Introducing $d_{3,n}$ with $n \geq 1$ such that
	\begin{equation}\label{2.17}
		d_{3,n}-qd_{3,n-1}=A_{3,n}, n\geq2.
	\end{equation}
	Hence, \eqref{2.16} can be written as
	\begin{equation}\label{2.18}
		b_{3,n+1}-tb_{3,n}-d_{3,n}=q(b_{3,n}-tb_{3,n-1}-d_{3,n-1}), n\geq2.
	\end{equation}
	Let  \begin{equation}\label{22.20}
		d_{3,0}=-\sum\limits_{n=1}^{\infty}\frac{A_{3,n}}{q^n}.
	\end{equation}
	Iterating \eqref{2.17} we have
	\begin{align}\nonumber
		d_{3,n}&=q^n(-\frac{A_{n+1}}{q^{n+1}}-\frac{A_{n+2}}{q^{n+2}}\cdots)\\\nonumber
		&=-\sum_{i=1}^{\infty}\frac{2t^{n+i}[D_{3,3}(n+i)^3+D_{3,2}(n+i)^2+D_{3,1}(n+i)+D_{3,0}]}{q^i}\\\nonumber
		&=-2t^n\sum_{i=1}^{\infty}(\frac{t}{q})^i[D_{3,3}(n+i)^3+D_{3,2}(n+i)^2+D_{3,1}(n+i)+D_{3,0}]\\
		&=:2t^n(D_{3,3}'n^3+D_{3,2}'n^2+D_{3,1}'n+D_{3,0}')\label{2.19} .
	\end{align}
	Now we solve
	\begin{equation}\label{2.20}
		\left\{
		\begin{array}{ll}
			4b_{3,1}+b_{3,2}&=A_{3,1},  \\
			b_{3,2}-tb_{3,1}&=d_{3,1}.
		\end{array}
		\right.
	\end{equation}
	Since the determinant of the coefficient matrix of the system \eqref{2.20} is non-zero, so  system  \eqref{2.20} has a unique solution for $b_{3,1}$. \\
	Iterating $b_{3,n+1}-tb_{3,n}-d_{3,n}=0$, we obtain
	\begin{align}
		b_{3,n}=tb_{3,n-1}+d_{3,n-1}&=t^{n-1}b_{3,1}+t^{n-2}d_{3,1}+t^{n-3}d_{3,2}+\cdots+d_{3,n-1}\nonumber\\
		&=t^{n-1}b_{3,1}+\sum_{i=1}^{n-1}(t^{n-1-i}d_{3,i})\nonumber\\
		&=t^{n-1}b_{3,1}+\sum_{i=1}^{n-1}[t^{n-1-i}2t^i(D_{3,3}'i^3+D_{3,2}'i^2+D_{3,1}'i+D_{3,0}')]\nonumber\\
		&=t^{n-1}b_{3,1}+2t^{n-1}\sum_{i=1}^{n-1}(D_{3,3}'i^3+D_{3,2}'i^2+D_{3,1}'i+D_{3,0}')\label{2.21}.
	\end{align}
	To simplify the expression \eqref{2.21}, we introduce the following lemma
	\begin{lemma}\label{th-2.1}  for any positive integer $k \in \mathrm{N}_{+}$, the sum expression $\displaystyle{}splaystyle{\sum\limits_{i=1}^ni^k}$  is a k+1 degree polynomial with respect to n, that is, exisits constant $C_0,C_1,\cdots,C_{k+1}$, such that
		$$\sum\limits_{i=1}^ni^k=C_{k+1}n^{k+1}+C_kn^k+\cdots+C_1n+C_0.$$

	\end{lemma}
	{\bf Proof.} By Stolz theorem,
	$$\lim\limits_{n \rightarrow \infty}\frac{\sum\limits_{i=1}^ni^k}{n^{k+1}}=\lim\limits_{n\rightarrow \infty}\frac{n^k}{n^{k+1}-(n-1)^{k+1}}=\lim\limits_{n\rightarrow \infty}\frac{n^k}{(k+1)n^{k}+\cdots+(-1)^k}=\frac{1}{k+1}.$$
	So there must exisits constant $C_0,C_1,\cdots,C_{k+1}$, such that
	$$\sum\limits_{i=1}^ni^k=C_{k+1}n^{k+1}+C_kn^k+\cdots+C_1n+C_0.$$
	The proof of Lemma \ref{th-2.1} is complete.\\
	
	By Lemma \ref{th-2.1},  \eqref{2.21} can be transformed into
	\begin{align}
		b_{3,n}&=t^{n-1}b_{3,1}+2t^{n-1}(D_{3,4}''n^4+D_{3,3}''n^3+D_{3,2}''n^2+D_{3,1}''n+D_{3,0}'')\nonumber\\
		&=:t^n(C_{3,4}n^4+C_{3,3}n^3+C_{3,2}n^2+C_{3,1}n+C_{3,0}). \label{2.22}
	\end{align}
	From the definition of $b_{3,n}$ given by \eqref{2.22}, we obtain
	$u_3(x)=\sum\limits_{n=1}^{\infty}b_{3,n}\sin nx$
	is absolutely and uniformly convergent on $[0,2\pi]$.\\
	
	Based on the expressions of $u_1$, $u_2$, $u_3$ obtained above, we summarize as follows.
	\begin{center}
		\begin{tabular}{|c|c|c|}
			\hline
			& $A_{m,n}$                                     & the Fourier  coefficients of $u_m$ \\
			\hline
			$m=1$ & --                                             & $2t^n$ \\
			$m=2$ & $2t^n(3n-1-\frac{2t^2}{1-t^2})$               & $-t^n(\frac{\sqrt{3}}{2}n^2+\frac{n}{3}+\frac{\sqrt{3}}{18})$ \\
			$m=3$ & $2t^n(D_{3,3}n^3+D_{3,2}n^2+D_{3,1}n+D_{3,0})$ & $t^n(C_{3,4}n^4+C_{3,3}n^3+C_{3,2}n^2+C_{3,1}n+C_{3,0})$ \\
			\hline
		\end{tabular}
	\end{center}
	
	\subsection{The expressions of $\mathbf{u_k}$ with $\mathbf{k> 3}$.}
	In this subsection, we give the expressions of $u_k$ with integers $k>3$.
	\begin{theorem}\label{th-3.1}
		Let $m$ be any positive integer. Assume that for any $k \in \mathrm{N}_{+}$ and $k\leq m-1$, one has
		\begin{equation}
			\begin{split}\label{3.1}
				u_k=&\frac{1+(-1)^k}{2}[\sum_{n=0}^{\infty}t^n(C_{k,2k-2}n^{2k-2}+C_{k,2k-3}n^{2k-3}+\cdots +C_{k,0})\cos nx]\\
				+&\frac{1-(-1)^k}{2}[\sum_{n=0}^{\infty}t^n(C_{k,2k-2}n^{2k-2}+C_{k,2k-3}n^{2k-3}+\cdots +C_{k,0})\sin nx],
			\end{split}
		\end{equation}
		where $C_{k,j}$ are some coefficients with $j=0,1,2,\cdots,2k-2$.
		Then when $k=m,$ the expression \eqref{3.1} is also ture.
	\end{theorem}
	
	{\bf Proof.} Let $p+q=m,$ where $p,q$ are positive integers.
	By \eqref{2.2}, one has
	\begin{equation}\label{3.2}
		u_0(x)u_m(x)=u_{m-1}^\prime(x)-\frac{1}{2}\sum_{p+q=m}u_p(x)u_q(x) .
	\end{equation}
	The following discusses the parity of $m$, mainly focusing on the Fourier coefficients of $u_pu_q$.
	\begin{enumerate}
		\item[Case 1.] $m$ is even.\\
		$(1)$ when both $p$ and $q$ are even, since $p,q\leq m-1$, we can assume that
		\begin{equation*}
			u_p(x)=\sum_{n=0}^{\infty}t^n(C_{p,2p-2}n^{2p-2}+C_{p,2p-3}n^{2p-3}+\cdots +C_{p,0})\cos nx,
		\end{equation*}
		\begin{equation*}
			u_q(x)=\sum_{n=0}^{\infty}t^n(C_{q,2q-2}n^{2q-2}+C_{q,2q-3}n^{2q-3}+\cdots +C_{q,0})\cos nx.
		\end{equation*}
		Denote
		\begin{align}
			a_{p,n}=t^n(C_{p,2p-2}n^{2p-2}+C_{p,2p-3}n^{2p-3}+\cdots +C_{p,0})\label{3.3}, \\
			a_{q,n}=t^n(C_{q,2q-2}n^{2q-2}+C_{q,2q-3}n^{2q-3}+\cdots +C_{q,0})\label{3.4}.
		\end{align}
		Then the coefficient of $\cos nx $ is:
		\begin{itemize}
			\item the coefficient of constant term is$\displaystyle{}splaystyle{\frac{a_{p,0}a_{q,0}}{2}+\sum\limits_{j=0}^{\infty}\frac{a_{p,j}a_{q,j}}{2}}; $
			\item the coefficient of $\cos nx$ is $\displaystyle{}splaystyle{\frac{1}{2}\sum\limits_{j=0}^{\infty}(a_{p,j}a_{q,j+n}+a_{p,j+n}a_{q,j})+\frac{1}{2}\sum\limits_{j=0}^{n}(a_{p,j}a_{q,n-j})}. $
		\end{itemize}
		Substituting \eqref{3.3} and \eqref{3.4} into the above expressions, we can obtain that the coefficient of $\cos nx$ is
		\begin{equation}
			\begin{split}\nonumber
				&\frac{1}{2}t^n\sum_{j=0}^{\infty}t^{2j}[(C_{p,2p-2}j^{2p-2}+C_{p,2p-3}j^{2p-3}+\cdots+C_{p,0})(C_{q,2q-2}(j+n)^{2q-2}+\cdots+C_{q,0})\\
				+&(C_{p,2p-2}(j+n)^{2p-2}+C_{p,2p-3}(j+n)^{2p-3}+\cdots+C_{p,0})(C_{q,2q-2}j^{2q-2}+C_{q,2q-3}j^{2q-3}+\cdots+C_{q,0})]\\
				+&\frac{1}{2}t^n\sum_{j=0}^{n}(C_{p,2p-2}j^{2p-2}+\cdots+C_{p,0})(C_{q,2q-2}(n-j)^{2q-2}+C_{q,2q-3}(n-j)^{2q-3}+\cdots+C_{q,0}).
			\end{split}
		\end{equation}
		$(2)$ when both $p$ and $q$ are odd, similarly, we can obtain that  the coefficient of $\cos nx$ is
		\begin{equation}
			\begin{split}\nonumber
				&\frac{1}{2}t^n\sum_{j=1}^{\infty}t^{2j}[(C_{p,2p-2}j^{2p-2}+C_{p,2p-3}j^{2p-3}+\cdots+C_{p,0})(C_{q,2q-2}(j+n)^{2q-2}+\cdots+C_{q,0})\\
				+&(C_{p,2p-2}(j+n)^{2p-2}+C_{p,2p-3}(j+n)^{2p-3}+\cdots+C_{p,0})(C_{q,2q-2}j^{2q-2}+C_{q,2q-3}j^{2q-3}+\cdots+C_{q,0})]\\
				-&\frac{1}{2}t^n\sum_{j=1}^{n}(C_{p,2p-2}j^{2p-2}+\cdots+C_{p,0})(C_{q,2q-2}(n-j)^{2q-2}+C_{q,2q-3}(n-j)^{2q-3}+\cdots+C_{q,0}).
			\end{split}
		\end{equation}
		Note that
		\begin{equation*}
			u_{m-1}'=\sum_{n=0}^{\infty}t^n(C_{m-1,2m-4}n^{2m-3}+\cdots+C_{m-1,0}n)\cos nx.
		\end{equation*}
		The coefficient of $\cos nx$ in $u_{m-1}^\prime(x)$ is $t^n(C_{m-1,2m-4}n^{2m-3}+\cdots+C_{m-1,0}n). $ By Lemma \ref{th-2.1}, the coefficient of $\cos nx$ in $\sum\limits_{p+q=m}u_{p}(x)u_{q}(x)$ is a polynomial of degree $2m-3$ with respect to $n$. Hence, we can assume that the coefficient of $\cos nx$ on the right-hand side of \eqref{3.2} is $t^n(D_{m,2m-3}n^{2m-3}+\cdots+D_{m,0}). $
		Let $A_{m,n}=2t^n(D_{m,2m-3}n^{2m-3}+\cdots+D_{m,0}),\ n\geq1.\ $
		The coefficient of $\cos nx$ on the left-hand side of \eqref{3.2} is
		$2a_{m,n}+\frac{a_{m,n-1}}{2}+\frac{a_{m,n+1}}{2},\ n\geq1.\ $ When $n=0,$ it is a constant term on the right-hand of \eqref{3.2}.\\ Denote$$A_{m,0}:=\frac{1}{2}\sum\limits_{p+q=m}(\sum_{j=0}^\infty a_{p,j}a_{q,j}+\sum_{j=1}^\infty b_{p,j}b_{q,j}+a_{p,0}a_{q,0}).\ $$
		Since $m$ is even, let $u_m(x)=\frac{a_{m,0}}{2}+\sum\limits_{n=1}^\infty a_{m,n}\cos nx.\ $ Then we have $$4a_{m,n}+a_{m,n-1}+a_{m,n+1}=A_{m,n},\quad n \geq 1.$$ Introduce $d_{m,n}$ satisfying
		\begin{equation}\label{3.5}
			\left\{
			\begin{array}{ll}
				a_{m,n+1}-ta_{m,n}-d_{m,n}=q(a_{m,n}-ta_{m,n-1}-d_{m,n-1}), \\
				d_{m,n}-qd_{m,n-1}=A_{m,n},\ n\geq1.\ \\
			\end{array}
			\right.
		\end{equation}
		Then by the second equation of \eqref{3.5}, we have
		\begin{align}
			d_{m,n}=qd_{m,n-1}+A_{m,n}&=q^2d_{m,n-2}+qA_{m,n-1}+A_{m,n}\nonumber\\
			&=q^nd_{m,0}+q^{n-1}A_{m,1}+q^{n-2}A_{m,2}+\cdots+qA_{m,n-1}+A_{m,n}\nonumber\\
			&=q^n[d_{m,0}+\frac{A_{m,1}}{q}+\frac{A_{m,2}}{q^2}+\cdots+\frac{A_{m,n}}{q^n}].\ \nonumber
		\end{align}
		Let $d_{m,0}=-\sum\limits_{k=1}^\infty \frac{A_{m,k}}{q^k}$. Then it follows that
		\begin{align}
			d_{m,n}&=-q^n(\frac{A_{m,n+1}}{q^{n+1}}+\cdots+\frac{A_{m,n+k}}{q^{n+k}}+\cdots)\nonumber\\
			&=-(\frac{A_{m,n+1}}{q}+\cdots+\frac{A_{m,n+k}}{q^k}+\cdots)\nonumber\\
			&=-\sum_{k=1}^\infty \frac{t^{n+k}[D_{m,2m-3}(n+k)^{2m-3}+\cdots+D_{m,0}]}{q^k}\nonumber\\
			&=-t^n\cdot\sum_{k=1}^\infty (\frac{t}{q})^k\cdot[D_{m,2m-3}(n+k)^{2m-3}+\cdots+D_{m,0}]\nonumber\\
			&=:t^n[D^\prime_{m,2m-3}n^{2m-3}+D^\prime_{m,2m-4}n^{2m-4}+\cdots+D^\prime_{m,0}].\ \nonumber
		\end{align}
		Solve
		\begin{equation*}
			\begin{cases}
				2a_{m,0}+a_{m,1}&=2A_{m,0}, \\
				a_{m,1}-ta_{m,0}&=d_{m,0}.
			\end{cases}
		\end{equation*}
		We obtain the unique solution $a_{m,0}$. Then by solving $a_{m,n}-ta_{m,n-1}=d_{m,n-1}$,
		we obtain
		\begin{align}
			a_{m,n}&=ta_{m,n-1}+d_{m,n-1}\nonumber\\
			&=t^na_{m,0}+t^{n-1}d_{m,0}+t^{n-2}d_{m,1}+\cdots+td_{m,n-2}+d_{m,n-1}\nonumber\\
			&=t^na_{m,0}+t^{n-1}D'_{m,2m-3}[1^{2m-3}+2^{2m-3}+\cdots+(n-1)^{2m-3}]\nonumber\\
			&~~+t^{n-2}D'_{m,2m-4}[1^{2m-4}+2^{2m-4}+\cdots+(n-1)^{2m-4}]+\cdots+D'_{m,0} n. \label{3.6}
		\end{align}
		Similarly, simplifying the last equation of \eqref{3.6}, we obtain by Lemma \ref{th-2.1}:
		\begin{equation*}
			\begin{split}
				a_{m,n}&=t^na_{m,0}+t^{n-1}\{D'_{m,2m-3}[1^{2m-3}+2^{2m-3}+\cdots+(n-1)^{2m-3}]\\
				&~~+D'_{m,2m-4}[1^{2m-4}+2^{2m-4}+\cdots+(n-1)^{2m-4}]+\cdots+D'_{m,0}\cdot n\}\\
				&=:t^n(C_{m,2m-2}n^{2m-2}+C_{m,2m-3}n^{2m-3}+\cdots+C_{m,1}n+C_{m,0}).
			\end{split}
		\end{equation*}
		The theorem holds true in this case.
		
		\item[Case 2.] $m$ is odd. In this case, $p$ and $q$ hold different parity.
		Without loss of generality, we assume that $p$ is odd and $q$ is even. Then, direct calculations yields
		\begin{itemize}
			\item the coefficient of $\sin x$ is $\frac{1}{2}[2a_{p,1}a_{q,0}+a_{p,2}a_{q,1}+a_{p,3}a_{q,2}+\cdots-(a_{p,1}a_{q,2}+a_{p,2}a_{q,3}+\cdots)]; $
			\item the coefficient of $\sin nx$ is
			\begin{equation}\label{3.7}
				\begin{split}
					&~~~~\frac{1}{2}(2a_{p,n}a_{q,0}+\sum_{j=1}^{n-1} a_{p,j}a_{q,n-j}+\sum_{j=n}^{\infty} a_{p,j+n}a_{q,j}-\sum_{j=1}^{\infty} a_{p,j}a_{q,j+n})\\
					=&\frac{1}{2}[2t^n(C_{p,2q-2}n^{2q-2}+\cdots+C_{p,0})C_{q,0}\\
					&+\sum_{j=1}^{n-1}t^j(C_{p,2p-2}j_{2p-2}+\cdots+C_{p,0})t^{n-j}(C_{q,2q-2}(n-j)^{2q-2}+\cdots+C_{q,0})\\
					&+\sum_{j=n}^{\infty}t^{n+j}(C_{p,2p-2}(n+j)^{2p-2}+\cdots+C_{p,0})t^j(C_{q,2q-2}j^{2q-2}+\cdots+C_{q,o})\\
					&-\sum_{j=1}^{\infty}t^j(C_{p,2p-2}j^{2p-2}+\cdots+C_{p,0})t^{n+j}(C_{q,2q-2}(n+j)^{2q-2}+\cdots+C_{q,0})],\ n\geq2.
				\end{split}
			\end{equation}
		\end{itemize}
		Moreover, we have $$u'_{m-1}=-\sum\limits_{n=0}^{\infty}t^n(C_{m-1,2m-4}n^{2m-3}+\cdots+C_{m-1,0}n)\sin nx.$$ By Lemma \ref{th-2.1}, the coefficient of $\sin nx$ in $\sum\limits_{p+q=m}u_{p}(x)u_{q}(x)$ is a polynomial of degree $2m-3$ with respect to $n$. Hence, we can assume that the coefficient of $\sin nx$ on the right-hand side of \eqref{3.7} is $t^n(D_{m,2m-3}n^{2m-3}+\cdots+D_{m,0}). $ Let $A_{m,n}=2t^n(D_{m,2m-3}n^{2m-3}+\cdots+D_{m,0}),\ n\geq1. $ The coefficient of $\sin nx$ on the left-hand side of \eqref{3.2} is
		\begin{itemize}
			\item the coefficient of $\sin x$ is $2b_{m,1}+\frac{b_{m,2}}{2}; $
			\item the coefficient of $\sin nx$ is $2b_{m,n}+\frac{b_{m,n-1}+b_{m,n+1}}{2}. $
		\end{itemize}
		Since $m$ is odd, we assume $u_m(x)=\sum\limits_{n=1}^{\infty} b_{m,n}\sin nx.$
		Comparing the coefficients on both sides of  \eqref{3.7} yields
		\begin{equation}\label{3.8}
			\begin{cases}
				4b_{m,1}+b_{m,2}=A_{m,1}, \\
				4b_{m,n}+b_{m,n-1}+b_{m,n+1}=A_{m,n},\ n\geq2.
			\end{cases}
		\end{equation}
		Introduce $d_{m,n}$ satisfying
		\begin{equation}\label{3.9}
			\begin{cases}
				b_{m,n+1}-tb_{m,n}-d_{m,n}=q[b_{m,n}-tb_{m,n-1}-d_{m,n-1}],\ \\
				d_{m,n}-qd_{m,n-1}=A_{m,n}.
			\end{cases}
		\end{equation}
		By \eqref{3.9}, we have
		\begin{equation}\nonumber
			\begin{split}
				d_{m,n}&=qd_{m,n-1}+A_{m,n}\\
				&=q^2d_{m,n-2}+qA_{m,n-1}+A_{m,n}\\
				&=q^{n-1}d_{m,1}+q^{n-2}A_{m,2}+\cdots+qA_{m,n-1}+A_{m,n}\\
				&=q^{n-1}\left(d_{m,1}+\frac{A_{m,2}}{q}+\cdots+\frac{A_{m,n}}{a^n}\right).
			\end{split}
		\end{equation}
		Let $d_{m,1}=-\sum\limits_{k=2}^{\infty} \frac{A_{m,k}}{q^{k-1}}$. Then it follows that
		\begin{equation}\nonumber
			\begin{split}
				d_{m,n}&=q^{n-1}\left(\frac{-A_{m,n+1}}{q^n}+\frac{-A_{m,n+2}}{q^{n+1}}+\cdots\right)\\
				&=-\sum_{k=1}^{\infty}\frac{2t^{n+k}[D_{m,2m-3}(n+k)^{2m-3}+\cdots+D_{m,0}]}{q^k}\\
				&=-2t^n\sum_{k=1}^{\infty}\left(\frac{t}{q}\right)^k[D_{m,2m-3}(n+k)^{2m-3}+\cdots+D_{m,0}]\\
				&=:t^n(D'_{m,2m-3}n^{2m-3}+D'_{m,2m-4}n^{2m-4}+\cdots+D'_{m,0}).
			\end{split}
		\end{equation}
		Solving
		\begin{equation}\nonumber
			\begin{cases}
				4b_{m,1}+b_{m,2}=A_{m,1},\ \\
				b_{m,2}-tb_{m,1}=d_{m,1},\
			\end{cases}
		\end{equation}
		we obtain the unique solution $b_{m,1}$. Then by solving $b_{m,2}=tb_{m,1}+d_{m,1},\ n\geq2$, we obtain
		\begin{equation}\nonumber
			\begin{split}
				b_{m,n}&=tb_{m,n-1}+d_{m,n-1}\\
				&=t^{n-1}b_{m,1}+t^{n-2}d_{m,1}+\cdots+td_{m,n-2}+d_{m,n-1}\\
				&=t^{n-1}b_{m,1}+t^{n-1}\{D'_{m,2m-3}[1^{2m-3}+2^{2m-3}\\
				&~~+\cdots+(n-1)^{2m-3}]+D'_{m,2m-4}[1^{2m-4}+2^{2m-4}+\cdots+(n-1)^{2m-4}]+\cdots+D'_{m,0}\cdot n\}\\
				&=t^n(C_{m,2m-2}n^{2m-2}+C_{m,2m-3}n^{2m-3}+\cdots+C_{m,1}n+C_{m,0}).
			\end{split}
		\end{equation}
	\end{enumerate}
	Therefore,  \eqref{3.1} holds for any $m\geq2$. The proof of Theorem \ref{th-3.1} is complete.
	
	\newpage
	\section{Proof of the main results}
	\subsection{First-order error estimates }
	In this subsection, we derive the first-order error equation and prove the existence of the periodic solutions to the error equation.\\
	
	Assume that $u^\v (x)=u_0(x)+\v v_1(\frac{x}{\v})$ satisfies \eqref{1.1}, then there hold
	\begin{equation*}
		-v_1^{\prime\prime}(\frac{x}{\v})+\v v_1(\frac{x}{\v})v_1^{\prime}(\frac{x}{\v})+
		u_0(x)v_1^{\prime}(\frac{x}{\v})+\v u_0^\prime(x)v_1(\frac{x}{\v})-
		u_0^{\prime\prime}(x)=0.
	\end{equation*}
	Let $\displaystyle{}{y=\frac{x}{\v}}$. It follows that
	\begin{equation*}
		-v_1^{\prime\prime}(y)+\v v_1(y)v_1^{\prime}(y)+u_0(\v y)v_1^{\prime}(y)+
		\v u_0^\prime(\v y)v_1(y)-u_0^{\prime\prime}(\v y)=0.
	\end{equation*}
	Note that $u_0(x)=-(2+\cos{x})$, it follows that
	\begin{equation}\label{4.1}
		-v_1^{\prime\prime}(y)+\v v_1(y)v_1^{\prime}(y)-(2+\cos(\v y))v_1^{\prime}(y)+
		\v \sin(\v y))v_1(y)-\v\cos(\v y)=0.
	\end{equation}
	For the convenience of presentation, we set $\v=\l^2$. Then \eqref{4.1} becomes
	\begin{equation}\label{4.2}
		-v_1^{\prime\prime}(y)+\l^2 v_1(y)v_1^{\prime}(y)-(2+\cos(\l^2 y))v_1^{\prime}(y)+
		\l^2 \sin(\l^2 y))v_1(y)-\l^2\cos(\l^2 y)=0.
	\end{equation}
	Then
		\begin{theorem}\label{th4.1}
		There exists $\l_{1}>0$ such that for any $0<| \l |<\l_{1}$, the equation \eqref{4.2}
		has a solution $v_{1}(y) \in C^2([0,\frac{2\pi}{\l^2}])$, satisfying $$v_{1}(0)=v_{1}(\frac{2\pi}{\l^2})$$
		and
		$$|v_{1}| \leq C,$$
		where $C>0$ is a constant independent of $\l$.
	\end{theorem}	
	
	To prove Theorem \ref{th4.1}, we first consider the initial problem of the equation \eqref{4.2}.
	\begin{lemma}
		Consider the Cauchy problem
		\begin{equation}\label{4.3}
			\left\{
			\begin{array}{ll}
				-v_1^{\prime\prime}(y)+\l^2 v_1(y)v_1^{\prime}(y)-(2+\cos(\l^2 y))v_1^{\prime}(y)+
				\l^2 \sin(\l^2 y))v_1(y)-\l^2\cos(\l^2 y)=0,\\
				v_1(0)=x_1,\\
				v_1^{\prime}(0)=y_1.
			\end{array}
			\right.
		\end{equation}
		where $\l \neq 0$, $x_1$ and $y_1$ are some given initial values. Then there exists $T_1>0$ such that the equation \eqref{4.3} has the unique solution $v_1(y) \in C^2[0,T_1)$.
	\end{lemma}
	
	{\bf Proof.} Note that the Cauchy problem \eqref{4.3} is equivalent to the following first-order
	ordinary differential equation
	\begin{equation}\label{4.4}
		\left\{
		\begin{array}{ll}
			-v_1^{\prime}(y)+\frac{\l^2}{2} v_1^2(y)-(2+\cos(\l^2 y))v_1(y)-
			\sin(\l^2 y)+y_1-\frac{\l^2}{2}x_1^2+3x_1=0,\\
			v_1(0)=x_1,\\
		\end{array}
		\right.
	\end{equation}
	Performing the classical Cauchy-Lipschitz for ODE, there exists $T_1>0$, such that \eqref{4.4}
	admits a unique local solution $v_1(y)\in C^1([0,T_1))$.$\hfill\Box$
	
	Next, we show that the existence interval $[0,T_1)$ can be large enough such that $T_1>\frac{2\pi}{\l^2}$.
	
	\begin{lemma}\label{La4.2}
		Assume that $0<|\l| <1$, and the initial data satisfy $x_1\leq-1$ and $\left| y_1\right| <1$, then the equation \eqref{4.4}
		admits a unique solution $v_1(y)\in C^1([0,\frac{2\pi}{\l^2}])$. Moreover, the solution $v_1(y)$
		is  bounded  uniformly with respect to $\l$.
	\end{lemma}
	
	{\bf Proof.} Let
	$$w_1(y)=e^{\int_0^y(2+\cos(\l^2t))dt}v_1(y)=
	e^{2y+\frac{1}{\l^2}\sin(\l^2y)}v_1(y).$$
	Then $w_1(y)$ satisfies
	\begin{equation}\label{4.5}
		\left\{
		\begin{array}{ll}
			w_1^{\prime}(y)=\frac{\l^2}{2}e^{-\big(2y+\frac{1}{\l^2}\sin(\l^2y)\big)}w_1^2(y)+
			e^{\big(2y+\frac{1}{\l^2}\sin(\l^2y)\big)}\big(y_1-\frac{\l^2}{2}x_1^2+3x_1-\sin(\l^2y)\big)\\
			w_1(0)=v_1(0)=x_1.
		\end{array}
		\right.
	\end{equation}
	The existence interval of the solution of the equation \eqref{4.4} is  same as that of the equation \eqref{4.5}, of which can be estimated by using upper and lower solutions method.
	
	Consider
	\begin{equation}\label{44.6}
		\left\{
		\begin{array}{ll}
			\overline{w}_1^{\prime}(y)=\frac{\l^2}{2}\overline{w}_1^2(y),\\
			\overline{w}_1(0)=x_1.
		\end{array}
		\right.
	\end{equation}
	Solve \eqref{44.6} to get
	$$\displaystyle{\overline{w}_1(y)=\frac{2x_1}{2-\l^2x_1y}},$$
	for $x_1<-1$ and $y \in [0,+\infty).$
	
	Consider
	\begin{equation}\label{44.7}
		\left\{
		\begin{array}{ll}
			\underline{w}_1^{\prime}(y)=(-x_1^2+3x_1-2)e^{\frac{5\pi}{\l^2}},\\
			\underline{w}_1(0)=x_1.
		\end{array}
		\right.
	\end{equation}
	Solve \eqref{44.7} to get
	$$\displaystyle{\underline{w}_1(y)=x_1+(-x_1^2+3x_1-2)e^{\frac{5\pi}{\l^2}}y}$$
	for $x_1<-1$ and $y \in [0,+\infty).$
	
	When $y \in[0,\frac{2\pi}{\l^2}]$, it follows that
	\begin{equation*}
		(-x_1^2+3x_1-2)e^{\frac{5\pi}{\l^2}}\leq
		\frac{\l^2}{2}e^{-\big(2y+\frac{1}{\l^2}\sin(\l^2y)\big)}w_1^2(y)+e^{\big(2y+\frac{1}{\l^2}\sin(\l^2y)\big)}\big(y_1-\frac{\l^2}{2}x_1^2+3x_1-\sin(\l^2y)\big)\leq
		\frac{\l^2}{2}w_1^2(y).
	\end{equation*}
	then we obtain $\underline{w}_1(y)\leq w_1(y)\leq \overline{w}_1(y)$ for any $y\in[0,\frac{2\pi}{\l^2}]$, which implies  $w_1(y)\in C^1([0,\frac{2\pi}{\l^2}]),v_1(y)\in C^1([0,\frac{2\pi}{\l^2}])$.
	
	Moreover, for any $0<|\l|<1$ and $y\in[0,\frac{2\pi}{\l^2}]$, we claim that $\displaystyle{}-x_1^2+3x_1-2\leq v_1(y)<0$. In fact
	\begin{equation*}
		v_1(y)=e^{-\big(2y+\frac{1}{\l^2}\sin(\l^2y)\big)}w_1(y)\leq e^{-\big(2y+\frac{1}{\l^2}\sin(\l^2y)\big)}\overline{w}_1=
		e^{-\big(2y+\frac{1}{\l^2}\sin(\l^2y)\big)}\frac{2x_1}{2-\l^2x_1y}<0.
	\end{equation*}
	and since
	\begin{equation*}
		v_1^{\prime}(y)=\frac{\l^2}{2} v_1^2(y)-(2+\cos(\l^2 y))v_1(y)-
		\sin(\l^2 y)+y_1-\frac{\l^2}{2}x_1^2+3x_1>-v_1(y)-x_1^2+3x_1-2.
	\end{equation*}
	There must be $v_1^{\prime}>0$, whenever $v_1<-x_1^2+3x_1-2$. Thus the claim is ture.\\
	The proof of Lemma \ref{La4.2} is complete.$\hfill\Box$
	
	Finally, we can complete the proof the Theorem \ref{th4.1}.\\
	{\bf Proof of Theorem 3.1.}
	We first give an equivalent condition of periodic solutions of \eqref{4.2}.\\
	Denote also the solution of \eqref{4.4} by $v_1(y,y_1)$. It follows from \eqref{4.4} that
	$$v_1(\frac{2\pi}{\l^2},y_1)=v_1(0,y_1)=x_1$$
	is equivalent to
	\begin{equation}\label{4.6}
		-v_1^{\prime}(\frac{2\pi}{\l^2},y_1)+\frac{\l^2}{2}v_1^2(\frac{2\pi}{\l^2},y_1)+y_1-\frac{\l^2}{2}x_1^2=0.
	\end{equation}
	Let $$\Phi_1(\l,y_1)=-v_1^{\prime}(\frac{2\pi}{\l^2},y_1)+\frac{\l^2}{2}v_1^2(\frac{2\pi}{\l^2},y_1)+y_1-\frac{\l^2}{2}x_1^2,$$
	where $(\l,y_1)\in U_1=\{(\l,y_1)\big|0< |\l|  <1, |y_1|<1\}$, $v_1$ is given by Lemma \ref{La4.2}.
	
	In order to determine the value of $\Phi_1$ at $\l= 0$, we consider
	\begin{equation}\label{44.9}
		\left\{
		\begin{array}{ll}
			-\tilde{v}_1^{\prime}(y,y_1)-3\tilde{v}_1(y,y_1)+3x_1+y_1=0,\\
			\tilde{v}_1(0,y_1)=x_1.
		\end{array}
		\right.		
	\end{equation}
	Solve \eqref{44.9}  to get
	$$\displaystyle{}\tilde{v}_1(y,y_1)=x_1+\frac{y_1}{3}(1-e^{-3y}).$$
	Let $r_1(y,y_1)=v_1(y,y_1)-\tilde{v}_1(y,y_1)$, then $r_1(y,y_1)$ satisfies
	\begin{equation}\label{44.10}
		\left\{
		\begin{array}{ll}
			-r_1^{\prime}(y,y_1)-3r_1(y,y_1)+\frac{\l^2}{2} v_1^2(y,y_1)+(1-\cos(\l^2 y))v_1(y,y_1)-
			\sin(\l^2 y)-\frac{\l^2}{2}x_1^2=0,\\
			r_1(0,y_1)=0.
		\end{array}
		\right.		
	\end{equation}
	Solve \eqref{44.10},  we obtain
	$$\displaystyle{} r_1(y,y_1)=e^{-3y}\Big(\int_{0}^{y}\big(\frac{\l^2}{2} v_1^2(t,y_1)+(1-\cos(\l^2 t))v_1(t,y_1)-
	\sin(\l^2 t)-\frac{\l^2}{2}x_1^2\big)e^{3t}dt\Big).$$
	In particular, for $\l >0,$
	\begin{align*}
		r_1(\frac{2\pi}{\l^2},y_1)&=e^{-\frac{6\pi}{\l^2}}\Big(\int_{0}^{\frac{2\pi}{\l^2}}\big(\frac{\l^2}{2} v_1^2(t,y_1)+(1-\cos(\l^2 t))v_1(t,y_1)-\sin(\l^2 t)-\frac{\l^2}{2}x_1^2\big)e^{3t}dt\Big)\\
		&=e^{-\frac{6\pi}{\l^2}}\Big(\int_{0}^{2\pi}\big(\frac{1}{2} v_1^2(\frac{t}{\l^2},y_1)+\frac{(1-\cos t)}{\l^2}v_1(\frac{t}{\l^2},y_1)-\frac{\sin t}{\l^2}-\frac{1}{2}x_1^2\big)e^{\frac{3t}{\l^2}}dt\Big).
	\end{align*}
	Since $v_1(y,y_1)$ is uniformly bounded with respect to $\l$ and $y_1$, we have
	$$\lim\limits_{\l \rightarrow 0}e^{-\frac{6\pi}{\l^2}}\int_{0}^{2\pi} Me^{\frac{3t}{\l^2}} dt =\lim\limits_{\l \rightarrow 0}\frac{M\l^2}{3}(1-e^{-\frac{6\pi}{\l^2}})=0,$$
	\begin{equation}\label{555}
		\lim\limits_{\l \rightarrow 0}e^{-\frac{6\pi}{\l^2}}\int_{0}^{2\pi}\frac{\sin t}{\l^2}e^{\frac{3t}{\l^2}} dt =\lim\limits_{\l \rightarrow 0}\frac{\l^2}{9+\l^4}(e^{-\frac{6\pi}{\l^2}}-1)=0,
	\end{equation}
	$$\lim\limits_{\l \rightarrow 0}e^{-\frac{6\pi}{\l^2}}\int_{0}^{2\pi}\frac{(1-\cos t)}{\l^2}e^{\frac{3t}{\l^2}} dt=\lim\limits_{\l \rightarrow 0}\frac{\l^4}{27+3\l^4}(1-e^{-\frac{6\pi}{\l^2}})=0,$$
	where $M$ is an arbitrary constant.\\
	Thus
	$$\lim\limits_{\l \rightarrow 0}v_1(\frac{2\pi}{\l^2},y_1)=\lim\limits_{\l \rightarrow 0}\tilde{v}_1(\frac{2\pi}{\l^2},y_1)=x_1+\frac{y_1}{3}.$$
	So according to \eqref{4.4}, we can define
	\begin{align*}
		\Phi_1(0,y_1)&=\lim\limits_{\l\rightarrow 0}\Phi_1(\l,y_1)\\		&=\lim\limits_{\l\rightarrow0}[-v_1^{\prime}(\frac{2\pi}{\l^2},y_1)+\frac{\l^2}{2}v_1^2(\frac{2\pi}{\l^2},y_1)+y_1-\frac{\l^2}{2}x_1^2]\\
		&=\lim\limits_{\l\rightarrow 0}[3v_1(\frac{2\pi}{\l^2},y_1)-3x_1]\\
		&=y_1
	\end{align*}
	for any $y_1 \in (-1,1).$ Then $\Phi_1(\l,y_1)$ is defined on  $\overline{U}_1=\{(\l,y_1)\big| |\l|  <1, |y_1|<1\}$. Now we verify its continuity on $\overline{U}_1.$ In particular, at the point $(0,0),$
	\begin{align*}
		\lim\limits_{y_1\rightarrow 0}\lim\limits_{\l \rightarrow 0}\Phi_1(\l,y_1)
		&=\lim\limits_{y_1\rightarrow 0}\lim\limits_{\l \rightarrow 0}[3v_1(\frac{2\pi}{\l^2},y_1)-3x_1]\\
		&=\lim\limits_{y_1\rightarrow 0}\lim\limits_{\l \rightarrow 0}[3\tilde{v}_1(\frac{2\pi}{\l^2},y_1)+3r_1(\frac{2\pi}{\l^2},y_1)-3x_1]\\
		&=\lim\limits_{y_1\rightarrow 0}\lim\limits_{\l \rightarrow 0}[y_1(1-e^{-\frac{6\pi}{\l^2}})]\\
		&=0.
	\end{align*}
	Observing that $\Phi_1(\l,y_1)$ converges to $y_1$ as $\l \rightarrow 0$ uniformly with respect to $y_1$, we have
	$$ \lim\limits_{\substack{y_1\rightarrow 0 \\ \l \rightarrow 0}}\Phi_1(\l,y_1)=\lim\limits_{y_1\rightarrow 0}\lim\limits_{\l \rightarrow 0}\Phi_1(\l,y_1)=0=\Phi_1(0,0). $$
	
	By the theorem on the differentiability of solutions of first-order ordinary differential equations with respect to parameters, $\frac{\partial \Phi_1}{\partial y_1}$ is continuous on $U_1$. In order to determine the value of $\frac{\partial \Phi_1}{\partial y_1}$ at $\l= 0$, we differentiate both sides of equation \eqref{4.4} with respect to $y_1$ yields
	\begin{equation}\label{4.11}
		\left\{
		\begin{array}{ll}
			-\frac{d}{d y}\frac{\partial v_1(y,y_1)}{\partial y_1}+{\l^2} v_1(y,y_1)\frac{\partial v_1(y,y_1)}{\partial y_1}-(2+\cos(\l^2 y))\frac{\partial v_1(y,y_1)}{\partial y_1}+1=0,\\
			\frac{\partial v_1}{\partial y_1}(0,y_1)=0.\\
		\end{array}
		\right.
	\end{equation}
	
	Next, we will use the method of upper and lower solutions to demonstrate that the solution $\frac{\partial v_1}{\partial y_1}(y,y_1)$ of equation \eqref{4.11} is uniformly bounded with respect to both $y$ and $y_1$. Since
	$$-x_1^2+3x_1-2 \leq v_1(y,y_1)<0$$
	for $y \in [0, \frac{2\pi}{\l^2}]$ and $y_1 \in (-1,1).$\\
	Consider
	\begin{equation}\label{4.12}
		\left\{
		\begin{array}{ll}
			\frac{d}{d y}\frac{\partial \overline{v}_1(y,y_1)}{\partial y_1}=-\frac{\partial \overline{v}_1(y,y_1)}{\partial y_1}+1,\\
			\frac{\partial \overline{v}_1}{\partial y_1}(0,y_1)=0.\\
		\end{array}
		\right.
	\end{equation}
	Solve \eqref{4.12} to get
	$$\frac{\partial \overline{v}_1}{\partial y_1}(y,y_1)=1-e^{-y}$$
	for $y_1 \in (-1,1)$ and $y \in [0,\infty).$

	Consider
	\begin{equation}\label{4.13}
		\left\{
		\begin{array}{ll}
			\frac{d}{d y}\frac{\partial \underline{v}_1(y,y_1)}{\partial y_1}=(-x_1^2+3x_1-5)\frac{\partial \underline{v}_1(y,y_1)}{\partial y_1}+1,\\
			\frac{\partial \underline{v}_1}{\partial y_1}(0,y_1)=0.\\
		\end{array}
		\right.
	\end{equation}
	Solve \eqref{4.13} to get
	$$\frac{\partial \underline{v}_1}{\partial y_1}(y,y_1)=\frac{1}{-x_1^2+3x_1-5}(e^{(-x_1^2+3x_1-5)y}-1)$$
	for $y_1 \in (-1,1)$ and $y \in [0,\infty).$

	Thus we have
	$$ 0 \leq \frac{\partial \underline{v}_1}{\partial y_1}(y,y_1) \leq \frac{\partial v_1}{\partial y_1}(y,y_1) \leq \frac{\partial \overline{v}_1}{\partial y_1}(y,y_1)\leq 1,$$
	for all $y \in [0, \frac{2\pi}{\l^2}]$ and $y_1 \in (-1,1)$, as desired.
	
	Because $v_1$ is continuously differentiable with respect to both $\l$ and $y_1$ in $U_1$, so
	\begin{align*}
		\frac{\partial r_1}{\partial y_1}(\frac{2\pi} {\l^2},y_1) &=e^{-\frac{6\pi}{\l^2}}\Big(\int_{0}^{\frac{2\pi}{\l^2}}\frac{\partial}{\partial y_1}\big(\frac{\l^2}{2} v_1^2(t,y_1)+(1-\cos(\l^2 t))v_1(t,y_1)-\sin(\l^2 t)-\frac{\l^2}{2}x_1^2\big)e^{3t}dt\Big)\\
		&=e^{-\frac{6\pi}{\l^2}}\Big(\int_{0}^{2\pi}\big(v_1^2(\frac{t}{\l^2},y_1)\frac{\partial v_1}{\partial y_1}(\frac{t}{\l^2},y_1)+\frac{(1-\cos t)}{\l^2}\frac{\partial v_1}{\partial y_1}(\frac{t}{\l^2},y_1)\big)e^{\frac{3t}{\l^2}}dt\Big).
	\end{align*}
	Given the uniform bounded of $\frac{\partial v_1}{\partial y_1}(\frac{2\pi}{\l^2},y_1)$ with respect to $\l$ and $y_1$, similar approach as in \eqref{555} yields
	$$\lim\limits_{\l \rightarrow 0}\frac{\partial r_1}{\partial y_1}(\frac{2\pi} {\l^2},y_1)=0$$
	uniformly in $y_1 \in (-1,1).$
	
	Define
	\begin{align*}
		\frac{\partial \Phi_1}{\partial y_1}(0,y_1)=\lim\limits_{\l \rightarrow 0} \frac{\partial \Phi_1}{\partial y_1}(\l^2,y_1)
		&=3\lim\limits_{\l \rightarrow 0} \frac{\partial v_1}{\partial y_1}(\frac{2\pi}{\l^2},y_1)\\
		&=3(\lim\limits_{\l \rightarrow 0} \frac{\partial \tilde{v}_1}{\partial y_1}(\frac{2\pi}{\l^2},y_1)+\lim\limits_{\l \rightarrow 0} \frac{\partial r_1}{\partial y_1}(\frac{2\pi}{\l^2},y_1))\\
		&=\lim\limits_{\l \rightarrow 0}(1-e^{-\frac{6\pi}{\l^2}})\\
		&=1
	\end{align*}
	for any $y_1 \in (-1,1).$

	Similar to the analysis of continuity above, $\frac{\partial \Phi_1(\l,y_1)}{\partial y_1}$ is continuous on $\overline{U}_1$.
	Since $$\Phi_1(0,0)=0,  \frac{\partial\Phi_1(0,0)}{\partial y_1}=1 \neq 0,$$
	by the implicit function theorem, there exists $0<\l_{1}<1$ such that for every $0<|\l|<\l_1$, there exists a unique $y_1$ satisfying $\Phi_1(\l,y_1(\l))=0$, that is \eqref{4.6} holds.
	The proof of Theorem \ref{th4.1} is complete.$\hfill\Box$

	\subsection{The (n + 1)-th order error estimates and proof of the main results} 	
	In this subsection, we derive the general error equation and prove the existence of the periodic solutions to the error equation and prove our main results.

	Assume that
	$$u^\v (x)=u_0(x)+\v u_1(x)+\v^2 u_2(x)+\dots+\v^n u_n(x)+\v^{n+1}v_{n+1}(\frac{x}{\v})$$
	satisfies \eqref{1.1}, then there holds
	\begin{equation*}
		\begin{aligned}
			&-v_{n+1}^{\prime\prime}(\frac{x}{\v})+\v^{n+1} v_{n+1}(\frac{x}{\v})v_{n+1}^{\prime}(\frac{x}{\v})+
			\big(u_0(x)+\v u_1(x)+\dots+\v^n u_n(x)\big)v_{n+1}^{\prime}(\frac{x}{\v})+\\
			&\big(\v u_0^\prime(x)+\v^2u_1^\prime(x)+\v^3 u_2^{\prime}(x)+\dots+\v^{n+1} u_n^{\prime}(x)\big) v_{n+1}(\frac{x}{\v})+\\
			&\sum_{k=1}^{n}{\v^k(u_k(x)u_n^\prime(x)+u_{k+1}(x)u_{n-1}^\prime(x)+\dots+u_n(x)u_k^\prime(x))}-\v u_n^{\prime\prime}(x)=0.
		\end{aligned}
	\end{equation*}
	Denote
	\begin{equation*}
		F(x)=\sum_{k=1}^{n}{\v^k(u_k(x)u_n^\prime(x)+u_{k+1}(x)u_{n-1}^\prime(x)+\dots+u_n(x)u_k^\prime(x))}-\v u_n^{\prime\prime}(x),
	\end{equation*}
	then $F(x)$ is a smooth and bounded function with period $2\pi$. \\
	
	Let $\displaystyle{}{y=\frac{x}{\v}}$. It follows that
	\begin{equation}\label{6.1}
		\begin{aligned}
			&-v_{n+1}^{\prime\prime}(y)+\v^{n+1} v_{n+1}(y)v_{n+1}^{\prime}(y)+\big(u_0(\v y)+\v u_1(\v y)+\dots+\v^n u_n(\v y)\big)v_{n+1}^{\prime}(y)\\
			&+\big(\v u_0^\prime(\v y)+\v^2u_1^\prime(\v y)+\v^3 u_2^{\prime}(\v y)+\dots+\v^{n+1} u_n^{\prime}(\v y)\big) v_{n+1}(y)
			+F(\v y)=0.
		\end{aligned}
	\end{equation}
	As in subsection 3.1, we set $\v=\l^2$. Then \eqref{6.1} becomes
	\begin{equation}\label{6.2}
		\begin{aligned}
			&-v_{n+1}^{\prime\prime}(y)+\l^{2n+2} v_{n+1}(y)v_{n+1}^{\prime}(y)+\big(u_0(\l^2 y)+\l^2 u_1(\l^2 y)+\dots+\l^{2n} u_n(\l^2 y)\big)v_{n+1}^{\prime}(y)\\
			&+\big(\l^2 u_0^\prime(\l^2 y)+\l^4u_1^\prime(\l^2 y)+\dots+\l^{2n+2} u_n^{\prime}(\l^2 y)\big) v_{n+1}(y)+F(\l^2 y)=0.
		\end{aligned}	
	\end{equation}
	It holds that
		\begin{theorem}\label{th6.1}
		There exists $\l_{n+1}>0$ such that for any $0< |\l|<\l_{n+1}$, the equation \eqref{6.2}
		has a solution $v_{n+1}(y) \in C^2([0,\frac{2\pi}{\l^2}])$, satisfying $$v_{n+1}(0)=v_{n+1}(\frac{2\pi}{\v})$$
		and
		$$|v_{n+1}| \leq C,$$
		where $C>0$ is a constant depending on $n$ but independent of $\v$.
	\end{theorem}	
	
	To prove Theorem \ref{th6.1}, we first consider the initial problem of equation \eqref{6.2}.
	\begin{lemma}
		Consider the Cauchy problem
		\begin{equation}\label{6.3}
			\left\{
			\begin{array}{l}
				\begin{aligned}
				&-v_{n+1}^{\prime\prime}(y)+\l^{2n+2} v_{n+1}(y)v_{n+1}^{\prime}(y)
				+\big(u_0(\l^2 y)+\l^2 u_1(\l^2 y)+\dots+\l^{2n} u_n(\l^2 y)\big)v_{n+1}^{\prime}(y)\\
				&+\big(\l^2 u_0^\prime(\l^2 y)+\l^4u_1^\prime(\l^2 y)+\dots+\l^{2n+2} u_n^{\prime}(\l^2 y)\big) v_{n+1}(y)+F(\l^2 y)=0,\\
					
				&v_{n+1}(0)=x_{n+1},\\
				&v_{n+1}^{\prime}(0)=y_{n+1}.
				\end{aligned}	\\
			\end{array}
			\right.
		\end{equation}
		where $\l \neq 0$, $x_{n+1}$ and $y_{n+1}$ are some given initial values. Then there exists $T_{n+1}>0$ such that  \eqref{6.3} has the unique solution $v_{n+1}(y) \in C^2([0,T_{n+1}))$.
	\end{lemma}
	
	{\bf Proof.} Note that the Cauchy problem \eqref{6.3} is equivalent to the following first-order
	ordinary differential equation
	\begin{equation}\label{6.4}
		\left\{
		\begin{array}{ll}
			\begin{aligned}
				&-v_{n+1}^{\prime}(y)+\frac{\l^{2n+2}}{2} v_{n+1}^2(y)+\big(u_0(\l^2 y)+\l^2 u_1(\l^2 y)+\l^4 u_2(\l^2 y)+\dots+\l^{2n} u_n(\l^2 y)\big)v_{n+1}(y)\\
			&	+F_1(\l^2 y)+y_{n+1}-\frac{\l^{2n+2}}{2}x_{n+1}^2-\big(u_0(0)+\l^2 u_1(0)+\l^4 u_2(0)+\dots+\l^{2n} u_n(0)\big)x_{n+1}=0,\\
				&v_{n+1}(0)=x_{n+1},\\
			\end{aligned}	\\
		\end{array}
		\right.
	\end{equation}
	where $\displaystyle{} F_1(\l^2 y)=\int_{0}^{\l ^2y}F(t)dt$.
	Performing the classical Cauchy-Lipschitz for ODE, there exists $T_{n+1}>0$ such that the equation \eqref{6.4}
	admits a unique local solution $v_{n+1}(y)\in C^1([0,T_{n+1}))$. $\hfill\Box$
	
	\begin{lemma}\label{La6.2}
		There exists $X_{n+1}>0$ and $\l_{n+1}^{\prime}>0$ such that for every $x_{n+1}<-X_{n+1},|y_{n+1}|<1,0<|\l|<\l_{n+1}^{\prime}$
		the equation \eqref{6.4}
		admits a unique solution $v_{n+1}(y)\in C^1([0,\frac{2\pi}{\l^2}])$. Moreover, the solution $v_{n+1}(y)$
		is  bounded  uniformly with respect to $\l$.
	\end{lemma}
	
	{\bf Proof.} By the boundedness and smoothness of $u_1,u_2\dots u_n$, and $u_0(y)=-(2+\cos y)$, there exists $0<\l_{n+1}^{\prime}<1$ and $X_{n+1}>0$, such that for any $y\in[0,\infty)$, if $0<|\l|<\l_{n+1}^{\prime}$ and $x_{n+1}<-X_{n+1}$, then
	$$-4<u_0(\l^2 y)+\l^2 u_1(\l^2 y)+\l^4 u_2(\l^2 y)+\dots+\l^{2n} u_n(\l^2 y)<-\frac{1}{2}$$
	and
	$$F_1(\l^2 y)+y_{n+1}-\frac{\l^{2n+2}}{2}x_{n+1}^2-\big(u_0(0)+\l^2 u_1(0)+\l^4 u_2(0)+\dots+\l^{2n} u_n(0)\big)x_{n+1}<0.$$
	Let $M_{n+1}>0$ such that
	$$-M_{n+1}<F_1(\l^2 y)+y_{n+1}-\frac{\l^{2n+2}}{2}x_{n+1}^2-\big(u_0(0)+\l^2 u_1(0)+\l^4 u_2(0)+\dots+\l^{2n} u_n(0)\big)x_{n+1}$$ for any $y\in[0,\infty)$.\\
	Let
	$$\displaystyle{} w_{n+1}(y)=e^{-\int_0^y\big(u_0(\l^2 t)+\l^2 u_1(\l^2 t)+\l^4 u_2(\l^2 y)+\dots+\l^{2n} u_n(\l^2 y)\big)dt}v_{n+1}(y),$$
	then $w_{n+1}(y)$ satisfies
	\begin{equation}\label{6.5}
		\left\{
		\begin{array}{l}
			\begin{aligned}
				w_{n+1}^{\prime}(y)=\frac{\l^{2n+2}}{2}e^{\int_0^y\big(u_0(\l^2 t)+\l^2 u_1(\l^2 t)+\dots+\l^{2n} u_n(\l^2 y)\big)dt}w_{n+1}^2(y)+\\
				e^{-\int_0^y\big(u_0+\l^2 u_1+\dots+\l^{2n} u_n\big)dt}
				\big(F_1(\l^2 y)+y_{n+1}-\frac{\l^{2n+2}}{2}x_{n+1}^2-\big(u_0(0)+\dots+\l^{2n} u_n(0)\big)x_{n+1}\big)\\
				w_{n+1}(0)=v_{n+1}(0)=x_{n+1}.
			\end{aligned}	\\
		\end{array}
		\right.
	\end{equation}
	The existence interval of the solution of the equation \eqref{6.4} is same as that of the equation \eqref{6.5}, of which can be estimated by upper and lower solutions method.
	
	More precisely, we consider
	\begin{equation}\label{65.6}
		\left\{
		\begin{array}{ll}
			\overline{w}_{n+1}^{\prime}(y)=\frac{\l^{2n+2}}{2}\overline{w}_{n+1}^2(y),\\
			\overline{w}_{n+1}(0)=x_{n+1}.
		\end{array}
		\right.
	\end{equation}
	Solve \eqref{65.6} to get
	$$\displaystyle{}\overline{w}_{n+1}(y)=\frac{2x_{n+1}}{2-\l^{2n+2}x_{n+1}y},$$
	for $x_{n+1}<-X_{n+1} $ and $ y\in[0,\infty).$
	
	Consider
	\begin{equation}\label{65.7}
		\left\{
		\begin{array}{ll}
			\underline{w}_{n+1}^{\prime}(y)=-M_{n+1}e^{\frac{8\pi}{\l^2}},\\
			\underline{w}_{n+1}(0)=x_{n+1}.
		\end{array}
		\right.
	\end{equation}
	Solving \eqref{65.7}, we obtain
	$$\displaystyle{}\underline{w}_{n+1}(y)=x_{n+1}-M_{n+1}e^{\frac{8\pi}{\l^2}}y,$$
	for $x_{n+1}<-X_{n+1} $ and $ y\in[0,\infty).$
	
	When $y \in[0,\frac{2\pi}{\l^2}]$, it follows that
	\begin{equation*}
		\begin{aligned}
			-M_{n+1}e^{\frac{8\pi}{\l^2}}\leq\frac{\l^{2n+2}}{2}e^{\int_0^y\big(u_0(\l^2 t)+\l^2 u_1(\l^2 t)+\dots+\l^{2n} u_n(\l^2 y)\big)dt}w_{n+1}^2(y)+
			e^{-\int_0^y\big(u_0+\dots+\l^{2n} u_n\big)dt}\\
			\big(F_1(\l^2 y)+y_{n+1}-\frac{\l^{2n+2}}{2}x_{n+1}^2-\big(u_0(0)+\dots+\l^{2n} u_n(0)\big)x_{n+1}\big)\leq\frac{\l^{2n+2}}{2}{w^2_{n+1}}(y).
		\end{aligned}
	\end{equation*}
	then we obtain $\underline{w}_{n+1}(y)\leq w_{n+1}(y)\leq \overline{w}_{n+1}(y)$ for any $y\in[0,\frac{2\pi}{\l^2}]$, which implies $w_{n+1}(y)\in C^1([0,\frac{2\pi}{\l^2}])$, $v_{n+1}(y)\in C^1([0,\frac{2\pi}{\l^2}])$.
	
	Moreover, for any $0<|\l|<\l_{n+1}^{\prime}$ and $y\in[0,\frac{2\pi}{\l^2}]$, we claim that $\displaystyle{}-2M_{n+1}\leq v_{n+1}(y)<0$. In fact,

	\begin{equation*}
		\begin{aligned}
			v_{n+1}(y)=e^{\int_0^y(u_0(\l^2 t)+\l^2 u_1(\l^2 t)+\dots+\l^{2n} u_n(\l^2 y))dt}w_{n+1}(y)\leq \\e^{\int_0^y(u_0(\l^2 t)+\dots+\l^{2n} u_n(\l^2 y))dt}\overline{w_{n+1}}(y)=
			e^{\int_0^y(u_0(\l^2 t)+\dots+\l^{2n} u_n(\l^2 y))dt}\frac{2x_{n+1}}{2-\l^{2n+2}x_{n+1}y}<0.
		\end{aligned}
	\end{equation*}
and since
	\begin{equation*}
		\begin{aligned}
			v_{n+1}^{\prime}(y)=\frac{\l^{2n+2}}{2} v_{n+1}^2(y)+\big(u_0(\l^2 y)+\l^2 u_1(\l^2 y)+\dots+\l^{2n} u_n(\l^2 y)\big)v_{n+1}(y)\\
			+F_1(\l^2 y)+y_{n+1}-\frac{\l^{2n+2}}{2}x_{n+1}^2-\big(u_0(0)+\l^2 u_1(0)+\dots+\l^{2n} u_n(0)\big)x_{n+1}>-\frac{1}{2}v_{n+1}(y)-M_{n+1}\\
		\end{aligned}
	\end{equation*}
	 there must be $v_{n+1}^{\prime}>0$, whenever $v_{n+1}<-2M_{n+1}$. Thus the claim is true.
	
	The proof of Lemma \ref{La6.2}  is complete.$\hfill\Box$
	
	Finally, we  complete the proof the Theorem \ref{th6.1}.\\
	{\bf Proof of Theorem \ref{th6.1}.}
	We first give an equivalent condition of periodic solutions.
	Denote also the solution of \eqref{6.4} by $v_{n+1}(y,y_{n+1})$.
	Setting $y=\frac{2\pi}{\l^2}$ in \eqref{6.4} and using the fact that $u_i (i=1,2,\cdots,n)$ are periodic function, it is observed that
	$$	v_{n+1}(\frac{2\pi}{\l^2},y_{n+1})=v_{n+1}(0,y_{n+1})$$
	is equivalent to
	\begin{equation}\label{6.6}
		-v_{n+1}^{\prime}(\frac{2\pi}{\l^2},y_{n+1})+\frac{\l^{2n+2}}{2}v_{n+1}^2(\frac{2\pi}{\l^2},y_{n+1})+y_{n+1}-\frac{\l^{2n+2}}{2}x_{n+1}^2=0.
	\end{equation}
	Let $$\Phi_{n+1}(\l,y_{n+1})=-v_{n+1}^{\prime}(\frac{2\pi}{\l^{2}},y_{n+1})+\frac{\l^{2n+2}}{2}v_{n+1}^2(\frac{2\pi}{\l^2},y_{n+1})+y_{n+1}-\frac{\l^{2n+2}}{2}x_{n+1}^2,$$
	where $(\l,y_{n+1})\in U_{n+1}=\{(\l,y_{n+1})\big|0< |\l|  <\l_{n+1}^{\prime}, |y_{n+1}|<1\}$, $v_{n+1}$ is given by Lemma \ref{La6.2} .
	
	In order to determine the value of $\Phi_{n+1}$ at $\l= 0$, we consider
	\begin{equation}\label{6.9}
		\left\{
		\begin{array}{ll}
			-\tilde{v}_{n+1}^{\prime}(y)-3\tilde{v}_{n+1}(y)+3x_{n+1}+y_{n+1}=0,\\
			\tilde{v}_{n+1}(0)=x_{n+1}.
		\end{array}
		\right.		
	\end{equation}
	Solve \eqref{6.9} to get
	$$\displaystyle{}\tilde{v}_{n+1}(y,y_{n+1})=x_{n+1}+\frac{y_{n+1}}{3}(1-e^{-3y}).$$
	Let $r_{n+1}(y,y_{n+1})=v_{n+1}(y,y_{n+1})-\tilde{v}_{n+1}(y,y_{n+1})$, then the $r_{n+1}(y,y_{n+1})$ satisfies the following equation
	\begin{equation}\label{6.10}
		\left\{
		\begin{array}{ll}
			-r_{n+1}^{\prime}(y,y_{n+1})-3r_{n+1}(y,y_{n+1})+\frac{\l^{2n+2}}{2} v_{n+1}^2(y,y_{n+1})+\big(3+u_0(\l^2 y)+\dots+\l^{2n} u_n(\l^2 y)\big)v_{n+1}\\
			+F_1(\l^2 y)-\frac{\l^{2n+2}}{2}x_{n+1}^2-\big(\l^2 u_1(0)+\l^4 u_2(0)+\dots+\l^{2n} u_n(0)\big)x_{n+1}=0,\\
			r_{n+1}(0)=0.
		\end{array}
		\right.		
	\end{equation}
	Solving \eqref{6.10}, we obtain
	\begin{equation*}
		\begin{aligned}
			r_{n+1}(y,y_{n+1})=e^{-3y}\Big(\int_{0}^{y}\big(\frac{\l^{2n+2}}{2} v_{n+1}^2(t,y_{n+1})+\big(3+u_0(\l^2 t)+\dots+\l^{2n} u_n(\l^2 t)\big)v_{n+1}(t,y_{n+1})\\
			+F_1(\l^2 t)-\frac{\l^{2n+2}}{2}x_{n+1}^2-\big(\l^2 u_1(0)+\l^4 u_2(0)+\dots+\l^{2n} u_n(0)\big)x_{n+1}\big)e^{3t}dt\Big).
		\end{aligned}
	\end{equation*}
	In particular,
	\begin{equation*}
		\begin{aligned}
			r_{n+1}(\frac{2\pi}{\l^2},y_{n+1})&=e^{-\frac{6\pi}{\l^2}}\Big(\int_{0}^{\frac{2\pi}{\l^2}}\big(\frac{\l^{2n+2}}{2} v_{n+1}^2(t,y_{n+1})+\big(3+u_0(\l^2 t)+\dots+\l^{2n} u_n(\l^2 t)\big)v_{n+1}(t,y_{n+1})\\
			&+F_1(\l^2 t)-\frac{\l^{2n+2}}{2}x_{n+1}^2-\big(\l^2 u_1(0)+\l^4 u_2(0)+\dots+\l^{2n} u_n(0)\big)x_{n+1}\big)e^{3t}dt\Big)\\
			&=e^{-\frac{6\pi}{\l^2}}\Big(\int_{0}^{2\pi}\big(\frac{\l^{2n}}{2} v_{n+1}^2(\frac{t}{\l^2},y_{n+1})+(\frac{1-\cos t}{\l^2}+u_1(t)+\cdots+\l^{2n-2} u_n(t))\cdot\\ &v_{n+1}(\frac{t}{\l^2},y_{n+1})
			+\frac{F_1(t)}{\l^2}-\frac{\l^{2n}}{2}x_{n+1}^2-\big(u_1(0)+\l^2 u_2(0)+\dots+\l^{2n-2} u_n(0)\big)x_{n+1}\big)e^{\frac{3t}{\l^2}}dt\Big).
		\end{aligned}
	\end{equation*}
	Since $v_{n+1}(y)$ is uniformly bounded with respect to $\l$ and $F_1(t)$ is smooth and periodic with $2\pi$.
	Then
	$$\lim\limits_{\l \rightarrow 0}e^{-\frac{6\pi}{\l^2}}\int_{0}^{2\pi} \frac{F_1(t)}{\l^2}e^{\frac{3t}{\l^2}} dt =\lim\limits_{\l \rightarrow 0}\frac{e^{-\frac{6\pi}{\l^2}}}{3}\int_{0}^{2\pi}F_1(t) d(e^{\frac{3t}{\l^2}})
	=\lim\limits_{\l \rightarrow 0}-\frac{e^{-\frac{6\pi}{\l^2}}}{3}\int_{0}^{2\pi}e^{\frac{3t}{\l^2}}F(t) dt=0.$$
	Thus
	$$\lim\limits_{\l \rightarrow 0}v_{n+1}(\frac{2\pi}{\l^2},y_{n+1})=\lim\limits_{\l \rightarrow 0}\tilde{v}_{n+1}(\frac{2\pi}{\l^2},y_{n+1})=x_{n+1}+\frac{y_{n+1}}{3}.$$
	So according to \eqref{6.4}, we can define
	\begin{align*}
		\Phi_{n+1}(0,y_{n+1})&=\lim\limits_{\l\rightarrow 0}\Phi_{n+1}(\l,y_{n+1})\\
		&=\lim\limits_{\l\rightarrow0}[-v_{n+1}^{\prime}(\frac{2\pi}{\l^2},y_{n+1})+\frac{\l^{2n+2}}{2}v_{n+1}^2(\frac{2\pi}{\l^2},y_{n+1})+y_{n+1}-\frac{\l^{2n+2}}{2}x_{n+1}^2]\\
		&=\lim\limits_{\l\rightarrow0}[(u_{0}(0)+\l^2u_1(0)+\cdots+\l^{2n}u_n(0))(x_{n+1}-v_{n+1}(\frac{2\pi}{\l^2},y_{n+1}))]\\
		&=\lim\limits_{\l\rightarrow 0}[3v_{n+1}(\frac{2\pi}{\l^{2}},y_{n+1})-3x_{n+1}]\\
		&=y_{n+1},
	\end{align*}
	for any $y_{n+1}\in (-1,1).$
	Then $\Phi_{n+1}(\l,y_{n+1})$ is defined on  $\overline U_{n+1}=\{(\l,y_{n+1})\big| |\l|  <\l_{n+1}^{\prime}, |y_{n+1}|<1\}$. Now we verify its continuity on $\overline{U}_{n+1}.$ In particular, at the point $(0,0),$
	\begin{align*}
		\lim\limits_{y_{n+1}\rightarrow 0}\lim\limits_{\l \rightarrow 0}\Phi_{n+1}(\l,y_{n+1})
		&=\lim\limits_{y_{n+1}\rightarrow 0}\lim\limits_{\l \rightarrow 0}[3v_{n+1}(\frac{2\pi}{\l^2},y_{n+1})-3x_{n+1}]\\
		&=\lim\limits_{y_{n+1}\rightarrow 0}\lim\limits_{\l \rightarrow 0}[3\tilde{v}_{n+1}(\frac{2\pi}{\l^2},y_{n+1})+3r_{n+1}(\frac{2\pi}{\l^2},y_{n+1})-3x_{n+1}]\\
		&=\lim\limits_{y_{n+1}\rightarrow 0}\lim\limits_{\l \rightarrow 0}[y_{n+1}(1-e^{-\frac{6\pi}{\l^2}})]\\
		&=0.
	\end{align*}
	Observing that $\Phi_{n+1}(\l,y_{n+1})$ converges to $y_{n+1}$ as $\l \rightarrow 0$ uniformly with respect to $y_{n+1}$, we have
	$$ \lim\limits_{\substack{y_{n+1}\rightarrow 0 \\ \l \rightarrow 0}}\Phi_{n+1}(\l,y_{n+1})=\lim\limits_{y_{n+1}\rightarrow 0}\lim\limits_{\l \rightarrow 0}\Phi_{n+1}(\l,y_{n+1})=0=\Phi_{n+1}(0,0). $$
	
	By the theorem on the differentiability of solutions of first-order ordinary differential equations with respect to parameters, $\frac{\partial \Phi_{n+1}}{\partial y_{n+1}}$ is continuous on $U_{n+1}$. In order to determine the value of $\frac{\partial \Phi_{n+1}}{\partial y_{n+1}}$ at $\l =0$, we differentiate both sides of equation \eqref{6.4} with respect to $y_{n+1}$ yields
	\begin{equation}\label{6.11}
		\left\{
		\begin{array}{ll}
			-\frac{d}{d y}\frac{\partial v_{n+1}(y,y_{n+1})}{\partial y_{n+1}}+{\l^{2n}} v_{n+1}(y,y_{n+1})\frac{\partial v_{n+1}(y,y_{n+1})}{\partial y_{n+1}}+(u_0(\l^2 y)+\cdots+\l^{2n} u_n(\l^2 y))\frac{\partial v_{n+1}(y,y_{n+1})}{\partial y_{n+1}}+1=0,\\
			\frac{\partial v_{n+1}}{\partial y_{n+1}}(0,y_{n+1})=0.\\
		\end{array}
		\right.
	\end{equation}
	
	Next, we will use the method of upper and lower solutions to demonstrate that the solution $\frac{\partial v_{n+1}}{\partial y_{n+1}}(y,y_{n+1})$ of equation \eqref{6.11} is uniformly bounded with respect to both $y$ and $y_{n+1}$. Since
	$$-2M_{n+1} \leq v_{n+1}(y,y_{n+1})<0, \qquad -4 < u_0(\l^2y)+\l^2u_1(\l^2y)+\cdots+\l^{2n}u_n(\l^2 y)<-\frac{1}{2}$$
	for $y \in [0, \frac{2\pi}{\l^2}]$ and $y_{n+1} \in (-1,1).$\\
	Consider
	\begin{equation}\label{6.12}
		\left\{
		\begin{array}{ll}
			\frac{d}{d y}\frac{\partial \overline{v}_{n+1}(y,y_{n+1})}{\partial y_{n+1}}=-\frac{1}{2}\frac{\partial \overline{v}_{n+1}(y,y_{n+1})}{\partial y_{n+1}}+1,\\
			\frac{\partial \overline{v}_{n+1}}{\partial y_{n+1}}(0,y_{n+1})=0.\\
		\end{array}
		\right.
	\end{equation}
	Solve \eqref{6.12} to get
	$$\frac{\partial \overline{v}_{n+1}}{\partial y_{n+1}}(y,y_{n+1})=2(1-e^{-2y})$$
	for $y_{n+1} \in (-1,1)$ and $y \in [0,\infty).$\\
	Consider
	\begin{equation}\label{6.13}
		\left\{
		\begin{array}{ll}
			\frac{d}{d y}\frac{\partial \underline{v}_{n+1}(y,y_{n+1})}{\partial y_{n+1}}=-(2M_{n+1}+4)\frac{\partial \underline{v}_{n+1}(y,y_{n+1})}{\partial y_{n+1}}+1,\\
			\frac{\partial \underline{v}_{n+1}}{\partial y_{n+1}}(0,y_{n+1})=0.\\
		\end{array}
		\right.
	\end{equation}
	Solve \eqref{6.13} to get
	$$\frac{\partial \underline{v}_{n+1}}{\partial y_{n+1}}(y,y_{n+1})=\frac{1}{2M_{n+1}+4}(1-e^{-\frac{y}{2M_{n+1}+4}})$$
	for $y_{n+1} \in (-1,1)$ and $y \in [0,\infty).$\\
	Thus we have
	$$ 0 \leq \frac{\partial \underline{v}_{n+1}}{\partial y_{n+1}}(y,y_{n+1}) \leq \frac{\partial v_{n+1}}{\partial y_{n+1}}(y,y_{n+1}) \leq \frac{\partial \overline{v}_{n+1}}{\partial y_{n+1}}(y,y_{n+1})\leq 2,$$
	for all $y \in [0, \frac{2\pi}{\l^2}]$ and $y_{n+1} \in (-1,1)$, as desired.\\
	
	Because $v_{n+1}$ is continuously differentiable with respect to both $\l$ and $y_{n+1}$ in $U_{n+1}$, so
	\begin{align*}
		\frac{\partial r_{n+1}}{\partial  y_{n+1}}(\frac{2\pi}{\l^2},y_{n+1})
		&=e^{-\frac{6\pi}{\l^2}}\Big(\int_{0}^{2\pi}\big({\l^{2n}} v_{n+1}(\frac{t}{\l^2},y_{n+1})\frac{\partial v_{n+1}}{\partial y_{n+1}}(\frac{t}{\l^2},y_{n+1})\\
		&+(\frac{1-\cos t}{\l^2}+u_1(t)+\cdots+\l^{2n-2}u_n(t))\frac{\partial v_{n+1}}{\partial y_{n+1}}(\frac{t}{\l^2},y_{n+1})
		\big)e^{\frac{3t}{\l^2}}dt\Big).
	\end{align*}
	Given the uniform bound of $\frac{\partial v_{n+1}}{\partial y_{n+1}}(\frac{2\pi}{\l^2},y_{n+1})$ with respect to $\l$ and $y_{n+1}$, similar qpproach as in \eqref{555} yields
	$$\lim\limits_{\l \rightarrow 0}\frac{\partial r_{n+1}}{\partial y_{n+1}}(\frac{2\pi} {\l^2},y_{n+1})=0$$
	uniformly in $y_{n+1} \in (-1,1).$\\
	
	We can define
	\begin{align*}
		\frac{\partial \Phi_{n+1}}{\partial y_{n+1}}(0,y_{n+1})=\lim\limits_{\l \rightarrow 0} \frac{\partial \Phi_{n+1}}{\partial y_{n+1}}(\l^2,y_{n+1})
		&=3\lim\limits_{\l \rightarrow 0} \frac{\partial v_{n+1}}{\partial y_{n+1}}(\frac{2\pi}{\l^2},y_{n+1})\\
		&=3(\lim\limits_{\l \rightarrow 0} \frac{\partial \tilde{v}_{n+1}}{\partial y_{n+1}}(\frac{2\pi}{\l^2},y_{n+1})+\lim\limits_{\l \rightarrow 0} \frac{\partial r_{n+1}}{\partial y_{n+1}}(\frac{2\pi}{\l^2},y_{n+1}))\\
		&=\lim\limits_{\l \rightarrow 0}(1-e^{-\frac{6\pi}{\l^2}})\\
		&=1
	\end{align*}
	for any $y_{n+1} \in (-1,1).$\\
	Similar to the analysis of continuity above, $\frac{\partial \Phi_{n+1}(\l,y_{n+1})}{\partial y_{n+1}}$ is continuous on $\overline{U}_{n+1}$. Since $$\Phi_{n+1}(0,0)=0,  \frac{\partial\Phi_{n+1}(0,0)}{\partial y_{n+1}}=1 \neq 0,$$
	by the implicit function theorem, there exists $0<\l_{n+1}<\l_{n+1}^{\prime}$ such that for every $0<|\l|<\l_{n+1}$, there exists a unique $y_{n+1}$ satisfying $\Phi_{n+1}(\l,y_{n+1}(\l))=0$, that is (\ref{6.6}) holds. \\
	The proof of Theorem \ref{th6.1} is complete.$\hfill\Box$
	
	The Proof of our main result Theorem \ref{th-11.2} then  follows from Theorem \ref{th6.1}.
	
	{\bf Acknowledgements.}
	This research was supported by the Interdisciplinary Project of Capital Normal University(Grant No. 2026JCYY04).

			\end{document}